\documentclass[amsmath,11pt]{article}

\usepackage{makeidx}         
\usepackage{graphicx}        
 \usepackage{color}                                
\usepackage{multicol}        
\usepackage[dvips]{epsfig}
\usepackage{graphics}
\usepackage{latexsym}
\usepackage{verbatim}
\usepackage{amssymb}
\usepackage{amsmath}
\usepackage{amsfonts}
\usepackage[tablesfirst,notablist]{endfloat}

\makeindex             
\newtheorem{theorem}{Theorem}[section]
\newtheorem{lemma}[theorem]{Lemma}
\newtheorem{remark}[theorem]{Remark}
\newtheorem{definition}[theorem]{Definition}

\def\real {{I \!\! R}}

\def\bigM{ {\mathcal {M}}}
\def\bigF{ {\mathcal {F}}}
\def\bigB{ {\mathcal {B}}}

\def\SS{ {\mathcal {S}}_{p,q}}
\def\RR{ {\mathcal {R}}_{p,q}}

\def\g{\hat g}

\newcommand{\rf}[1]{(\ref{#1})}

\def\la{\lambda}

\def\g{\frac{p^2+q^2-1} 2}
\def\e{\epsilon}
\def\a{\alpha}
\def\d{\delta}

\def\ff{\widehat f}

\def\({\left(}
\def\){\right)}

\def\be#1\ee{\begin{equation}#1\end{equation}}
\def\ba#1#2\ea{\begin{array}#1#2\end{array}}
\def\bgr#1\egr{{\allowdisplaybreaks\begin{gather}#1\end{gather}}}
\def\bma#1\ema{{\allowdisplaybreaks\begin{align}#1\end{align}}}

\def\oplem#1{\begin{lemma}\, {\rm #1}\, \it }
\def\cllem{\end{lemma}\rm \par }
\def\opthm#1{\begin{theorem}\, {\rm #1}\, \it }
\def\clthm{\end{theorem}\rm \par }

\def\R{\mathbb R}

\newcommand{\fer}[1]{(\ref{#1})}
\newcommand{\bq}{\begin{equation}}
\newcommand{\eq}{\end{equation}}
\def\bqa{\begin{eqnarray}}
\def\eqa{\end{eqnarray}}
\def\bd{\begin{displaymath}}
\def\ed{\end{displaymath}}

\begin{document}

\title{Self--similarity and power--like tails\\ in nonconservative kinetic models}
\author{Lorenzo Pareschi\thanks{Department of Mathematics and Center for Modelling Computing and Statistics (CMCS),
University of Ferrara, Via Machiavelli 35 I-44100 Ferrara, Italy.
\texttt{lorenzo.pareschi@unife.it}}\and Giuseppe
Toscani\thanks{Department of Mathematics, University of Pavia, via
Ferrata 1, 27100 Pavia, Italy.
\texttt{giuseppe.toscani@unipv.it}}}
%
%
\date{January 11, 2006}

\maketitle

\begin{abstract}
In this paper, we discuss the large--time behavior of solution of
a simple kinetic model of Boltzmann--Maxwell type, such that the
temperature is time decreasing and/or time increasing. We show
that, under the combined effects of the nonlinearity and of the
time--monotonicity of the temperature, the kinetic model has non
trivial quasi-stationary states with power law tails.  In order to
do this we consider a suitable asymptotic limit of the model
yielding a Fokker-Planck equation for the distribution. The same
idea is applied to investigate the large--time behavior of an
elementary kinetic model of economy involving both exchanges
between agents and increasing and/or decreasing of the mean
wealth. In this last case, the large--time behavior of the
solution shows a Pareto power law tail. Numerical results confirm
the previous analysis.
\end{abstract}

{\bf Keywords.} Granular gases, overpopulated tails, Boltzmann equation, wealth and income
distributions, Pareto distribution.

\tableofcontents

\section{Introduction}
\label{sec:1}
  A well--known phenomenon in the large--time behavior of the Boltzmann equation with
dissipative interactions  is the formation of overpopulated tails
\cite{BK, EB1, EB2}. Exact results on the behavior of these tails
have been obtained for simplified models, in particular for a gas
of inelastic Maxwell particles. Our goal here is to show that, at
least for some simplified kinetic model,  the formation of
overpopulated tails is not only a behavior typical of systems
where there is dissipation of the temperature (cooling), but more
generally is a consequence of the fact that the temperature is not
conserved.
One can indeed conjecture that the formation of overpopulated
tails in a kinetic model depends on the breaking of energy
conservation.
In kinetic theory of rarefied gases, formation of overpopulated
tails has been first observed for inelastic Maxwell models
\cite{EB1, EB2}. Inelastic Maxwell models share with elastic
Maxwell molecules the property that the collision rate in the
Boltzmann equation is independent of the relative velocity of the
colliding pair. These models are of interest for granular fluids
in spatially homogeneous states because of the mathematical
simplifications resulting from a velocity independent collision
rate. Among others properties, the inelastic Maxwell models
exhibit similarity solutions, which represent the intermediate
asymptotic of a wide class of initial conditions \cite{BCT}.
Recently, the study of a dissipative kinetic model obtained by
generalizing the classical model known as Kac caricature of a
Maxwell gas \cite{PT}, led to new ideas on the mechanism of the
formation of tails. Indeed, in \cite{PT} connections between the
cooling problem for the dissipative model and the classical
central limit theorem for stable laws of probability theory were
found. A second point in favor of our conjecture on tails
formation comes out from some recent applications to economy of
one--dimensional kinetic models of Maxwell type \cite{Sl, P, CPT}.
The main physical law here is that a strong economy produces
growth of the mean wealth (which of course is the opposite
phenomenon to the dissipation). Nevertheless, the kinetic model
led to an immediate explanation of the formation of Pareto tails
\cite{Par}. Having this in mind, in the next Section we study a
one--dimensional Boltzmann--like equation which is able to
describe both dissipation and production of energy. This model has
been recently considered in \cite{BBLR} with the aim of recovering
exact self-similar solutions. The analysis of \cite{BBLR}, based
on the possibility to use Fourier transform techniques to
investigate properties of the self-similar profiles, shows that in
many cases there is evidence of algebraic decay of the velocity
distributions. On the other hand, except in particular cases, no
exact results can be achieved. To obtain a almost complete
description of the large time behavior of the solution, we resort
to a different approach. After a brief description of the model,
in Section 2 we introduce a suitable asymptotic analysis, which
reduces the Boltzmann equation to a Fokker--Planck like equation
which has an explicitly computable stationary state with
power--like tails. In Section 3, we show how similar ideas can be
fruitfully applied to describe the large--time behavior of some
elementary kinetic models of an open economy. Here, the underlying
Fokker--Planck equation takes the form of a similar one introduced
recently in \cite{BM, CPT}. The rest of the paper is devoted to
the proof of mathematical details. Numerical experiments on the
Boltzmann models can be found at the end of the paper.

\section{Kinetic models and Fokker-Planck asymptotics}

 In this section we will study the large--time behavior of solutions to one--dimensional
 kinetic models of Maxwell-Boltzmann type, where the binary interaction between particles obey to
 the law
 \be\label{coll}
 v^* = pv+qw, \quad w^* = qv +pw ; \quad p>q>0 .
 \ee
The positive constants $p$ and $q$ represent the interacting
parameters, namely the portion of the
 pre--collisional velocities $(v,w)$ which generate the post--collisional ones $(v^*,w^*)$.
 As it will be clear after Subsection \ref{graz2}, the choice
 $p>q$ is natural in mimicking economic interactions, so that we
 will assume it even in molecular dynamics. As a matter of fact,
the mixing parameters $p$ and $q$ can be exchanged, which
corresponds to the exchange of post--collision velocities, without
any change in the global collision evolution.

\subsection{Nonconservative kinetic models}\label{graz1}
 Let  $f(v,t)$ denote the distribution of particles with velocity $v \in \real$ at time $t \ge
0$. The kinetic model can be easily derived by standard methods of
kinetic theory, considering that the change in time of $f(v,t)$
depends on a balance between the gain and loss of particles with
velocity $v$ due to binary collisions. This leads to the following
integro-differential equation of Boltzmann type \cite{BBLR},
 \be \frac{\partial f}{\partial t} =
\int_{\real}\left( \frac 1J f(v_*)
 f(w_*) - f(v) f(w)\right) d
 w \label{eq:boltz}
  \ee
  where $(w_*,w_*)$ are the
pre-collisional velocities that generate the couple $(v,w)$ after
the interaction. In (\ref{eq:boltz}) $J = p^2-q^2$ is the Jacobian
of the transformation of $(v,w)$ into $(v^*,w^*)$. Note that,
since we fixed  $p>q$, the Jacobian $J$ is positive and that the
unique situation corresponding to $J=1$ is obtained taking $p=1$
and $q=0$ for which the collision operator vanishes.

The kinetic equation \fer{eq:boltz} is the analogous of the
Boltzmann equation for Maxwell molecules \cite{Bob, CIP}, where
the collision frequency is assumed to be constant. Also, it
presents several similarities with the one-dimensional Kac model
\cite{Kac, McK}. It is well-known to people working in kinetic
theory that this simplification allows for a better understanding
of the qualitative behavior of the solutions.

Without loss of generality, we can fix the initial density to
satisfy
 \be\label{norm1}
 \int_{\real} f_0(v)\, dv =1 \, ; \quad    \int_{\real} v f_0(v)\, dv =  0 \,
 \quad  \int_{\real} v^2 f_0(v)\, dv = 1.
 \ee
To avoid the presence of the Jacobian, and to study approximation
to the collision operator
  it is extremely convenient to write equation \fer{eq:boltz} in weak form.
 It corresponds to consider, for all smooth functions $\phi(v)$,
 the equation
 \be\label{weak boltz}
\frac d{dt}\int_{\real} \phi(v)f(v,t)\,dv  =
 \int_{\real^2} f(v)
f(w) ( \phi(v^*)-\phi(v)) dv d w .
 \ee
One can alternatively use the symmetric form
 \begin{eqnarray} \label{weak boltz2}\nonumber
\frac{d}{dt} \int_{\real} f(v) \phi(v)\,dv &=& \frac12 \int_{\real^2}  f(v) f(w) \\
\\[-.25cm]
\nonumber &&( \phi(v^*)+\phi(w^*)-\phi(v)-\phi(w)) dv \, dw\, .
  \end{eqnarray}
  A remarkable fact is that
 equations \fer{weak boltz} and \fer{weak boltz2} can be studied
  for all values of the mixing parameters $p$ and $q$, including
  the case $p=q$, which could not be considered in equation
  \fer{eq:boltz}.

Choosing $\phi(v) = v$, (respectively $\phi(v) = v^2$) shows that
 \be
 m(t) = \int_{\real} v f(v,t) \,dv = m(0)\exp\left\{ (p+q-1)t\right\}.
 \ee
 Hence, since the initial density $f_0$ satisfies \fer{norm1}, $m(0) =0$ and  $m(t)
 =0$ for all $t>0$.
 Consequently,
 \be
 E(t) = \int_{\real} v^2 f(v,t)  \,dv = \exp\left\{ (p^2+q^2-1)t\right\}.
 \ee
Higher order moments can be evaluated recursively, remarking that
the integrals $\int v^nf(v,t)$ obey a closed hierarchy of
equations \cite{BK}.

Note that the second moment of the solution is not conserved,
unless the collision parameters satisfy
 \[
 p^2 + q^2 = 1.
 \]
If this is not the case, the energy can grow to infinity or
decrease to zero, depending on the sign of $p^2+q^2-1$. In both
cases, however,  stationary solutions of finite energy do not
exist, and the large--time behavior of the system can at best be
described by self-similar solutions. The standard way to look for
self--similarity is to scale the solution according to the role
 \be\label{scala}
  g(v,t) = \sqrt{ E(t)}f\left( v\sqrt{ E(t)}, t \right) .
  \ee
This scaling implies that $\int v^2 g(v,t)= 1$ for all $t \ge 0$.
Elementary computations show that $g= g(v,t)$ satisfies
 the equation
  \be \frac{\partial g}{\partial t}  -\frac 12\left(p^2 +q^2 -1\right)\frac{\partial}{\partial v}\left(vg\right)=
\int_{\real}\left( \frac 1J g(v_*)
 g(w_*) - g(v) g(w)\right) dw \label{boltz1} .
  \ee
In weak form, equation \fer{boltz1} reads
 \[\nonumber
\frac d{dt}\int_{\real} \phi(v)g(v,t)\,dv  -\frac 12\left(p^2 +q^2
-1\right)\int_{\real}\phi(v)\frac{\partial}{\partial v}\left(vg\right)\, dv =
 \]
 \be\label{weak boltz1}
 \int_{\real^2} g(v)
g(w) ( \phi(v^*)-\phi(v)) dv d w .
 \ee
Assuming that $\phi$ vanishes at infinity, we can integrate by
parts the second integral on the right--hand side of \fer{weak
boltz1} to obtain
 \[\nonumber
\frac d{dt}\int_{\real} \phi(v)g(v,t)\,dv +\frac 12\left(p^2 +q^2
-1\right)\int_{\real}\phi'(v)vg(v)\, dv =
 \]
  \be\label{weak boltz2b}
 = \int_{\real^2} g(v)
g(w) ( \phi(v^*)-\phi(v)) dv d w .
 \ee
By the collision rule (\ref{coll}),
$$
v^*   - v =   (p-1)v +qw .
$$
Let us use a second order Taylor expansion of $\phi(v^*)$ around
$v$
 $$
\phi(v^*) - \phi(v) = \left((p-1)v +qw\right) \phi'(v) + \frac 12\left((p-1)v +qw\right)^2 \phi''(\tilde v),
 $$
where, for some $0 \le \theta \le 1$
 \[
 \tilde v = \theta v^* +(1 -\theta)v .
  \]
Inserting this expansion in the collision operator,  we obtain the
equality
  \[
 \int_{\real^2} g(v)
g(w) ( \phi(v^*)-\phi(v)) dv d w =\int_{\real^2} g(v) g(w) \left((p-1)v +qw\right) \phi'(v)dv d w
+
 \]
  \be\label{taylor}
\frac 12\int_{\real^2} g(v) g(w)\left((p-1)v +qw\right)^2 \phi''(v) dv d w + R(p,q),
 \ee
where
 \be\label{resto}
 R(p, q) = \frac 1{2}\int_{\real^2}\left((p-1)v +qw\right)^2  \left(
\phi''(\tilde v)- \phi''( v)\right)
 g(v)g(w)  dv\,d w .
 \ee
  Recalling that $g(v,t)$ satisfies \fer{norm1}, we can simplify into \fer{taylor} to obtain
  \[
 \int_{\real^2} g(v)
g(w) ( \phi(v^*)-\phi(v)) dv d w = (p-1)\int_{\real} vg(v)\phi'(v)dv +
 \]
  \be\label{taylo1}
\frac 12\int_{\real} g(v)\left((p-1)^2v^2 +q^2\right)\phi''(v) dv + R(p,q).
 \ee
 Substituting \fer{taylo1} into \fer{weak boltz2b}, and grouping similar terms, we conclude that
 $g(v,t)$ satisfies
 \[
\frac d{dt}\int_{\real} \phi(v)g(v,t)\,dv +\frac 12\left((p-1)^2 +
q^2\right)\int_{\real}\phi'(v)vg(v)\, dv =
 \]
 \be
\frac 12\int_{\real} g(v)\left((p-1)^2v^2 +q^2\right)\phi''(v) dv + R(p,q).
 \ee
Hence, if we set
 \be\label{resc}
 \tau = q^2 t, \quad h(v,\tau) = g(v,t),
 \ee
 which implies $g_0(v) = h_0(v)$, $h(v,\tau)$ satisfies
  \[
\frac d{d\tau}\int_{\real} \phi(v)h(v,\tau)\,dv +\frac 12\left(\left(\frac{p-1}q\right)^2 +
1\right)\int_{\real}\phi'(v)vh(v)\, dv =
 \]
 \be\label{quasiFP}
\frac 12\int_{\real} h(v)\left(\left(\frac{p-1}q\right)^2v^2 +1\right)\phi''(v) dv + \frac 1{q^2}R(p,q).
 \ee
 Suppose now that the remainder in \fer{quasiFP} is small
 for small values of the parameter $q$. Then equation \fer{quasiFP} gives the behavior of
 $g(v,t)$ for large values of time.
Moreover, taking $p=p(q)$ such that, for a given constant $\lambda$
  \be\label{con}
 \lim_{q\to 0} \frac{p(q)-1}q = \lambda,
  \ee
  equation \fer{quasiFP} is well--approximated by the equation (in weak form)
 \[
\frac d{d\tau}\int_{\real} \phi(v)h(v,\tau)\,dv +\frac 12\left(\lambda^2 +
1\right)\int_{\real}\phi'(v)v h(v)\, dv =
 \]
 \be\label{FP}
 \frac 12\int_{\real} h(v)\left(\lambda^2v^2 +1\right)\phi''(v) dv.
 \ee
 Equation \fer{FP} is nothing but the weak form of the Fokker-Planck equation
 \be\label{FP1}
 \frac{\partial h}{\partial \tau} = \frac 12\left( \frac{\partial^2 }
 {\partial v^2}\left( (1 + \lambda^2v^2)h\right) + \left(1 + \lambda^2 \right)
 \frac{\partial }{\partial v}\left(v
 h\right)\right),
 \ee
which has a unique stationary state of unit mass, given by
 \be\label{staz}
  M_\lambda(v) = c_\lambda \left( \frac 1{1 + \lambda^2v^2}\right)^{\frac 32 +
  \frac1{2\lambda^2}},
  \ee
where
 \be\label{cl}
c_\lambda=\frac{|\lambda|}{\sqrt{\pi}}\frac{\displaystyle\Gamma
\left(\frac{3\lambda^2+1}{2\lambda^2}\right)}{\displaystyle\Gamma
\left(\frac{1+2\lambda^2}{2\lambda^2}\right)}.
\ee

  \begin{remark} The derivation of the Fokker-Planck equation
 \fer{FP1} presented in this section is largely formal. The main
 objective here was to show that there are regimes of the mixing
 parameters for which we can expect formation of self-similar
 solutions to the kinetic model with overpopulated tails. We postpone the detailed
 proof and the mathematical technicalities to the second part of the
 paper.
\end{remark}

\begin{remark}
The conservative case $p^2 + q^2 = 1$ can be treated likewise. In
this case one is forced to choose $p = \sqrt{1-q^2}$, which gives
$\lambda =0$ as unique possible value.  In the limit one then
obtains the linear Fokker--Planck equation
 \be\label{FP2}
 \frac{\partial h}{\partial \tau} = \frac 12\left( \frac{\partial^2 h}
 {\partial v^2} +
 \frac{\partial }{\partial v}\left(v
 h\right)\right).
 \ee
 Note that in this case the stationary solution $M(v)$ is the Maxwell density
 \be
 M(v) = \frac 1{\sqrt{2\pi}}e^{-v^2/2},
 \label{eq:max}
 \ee
 for all $q < 1/\sqrt{2}$.
 On the contrary, the non conservative cases are characterized by a $\lambda$ different from
 zero, which produces a stationary state with overpopulated tails.
 Note that from (\ref{cl}) we have $c_\lambda \to 1/\sqrt{2\pi}$
 as $\lambda\to 0$ and thus $M(v)=\lim_{\lambda\to 0}
 M_{\lambda}(v)$.
\end{remark}

\begin{remark}\label{poss}
 The possibility to pass to the limit in \fer{con}, with $\la >0$, is restricted to the cases $p^2+q^2 <1$
 and $p^2+q^2 >1$, but $p>1$. In the  case $p^2+q^2 >1$, $p<1$, it holds
  \[
 0 < \frac{(1-p)^2}{q^2} < \frac{1-p}{1+p},
  \]
 which forces $\la$ towards zero as $p \to 1$. This case, as the conservative one,
gives in the limit the linear Fokker--Planck equation. Hence, formation of tails is expected in case of dissipation
of energy, as well in case of production of energy, but only when the mixing parameter $p>1$.
 \end{remark}

\begin{remark}\label{granu}
In addition to the conservative case, a second one deserves to be mentioned. If $p = 1-q$, the kinetic models is
nothing but the model for granular dissipative collisions introduced and studied in \cite{McNY, BK, BMP} as a
one--dimensional caricature of the Maxwell--Boltzmann equation \cite{BCG, BC4}. In this case $\lambda = -1$, and
the stationary state is
 \be
 \label{eq:std1}
 M_1(v) = \frac 2\pi\left( \frac 1{1 + v^2}\right)^2.
 \ee
This solution solves the kinetic equation \fer{weak boltz1}, for any value of the parameter $q<1/2$.
\end{remark}

\begin{remark}
The asymptotic procedure considered in this section is the analogue of the so-called grazing collision limit of the
Boltzmann equation \cite{Vi, vil2}, which relies in concentrating the rate functions on collisions which are
grazing, so leaving the collisional velocities unchanged. It is well--known that in this (conservative) case, while
the Boltzmann equation changes into the Landau--Fokker--Planck equation, the stationary distribution remains of
Maxwellian type.
\end{remark}

\subsection{Pareto tails in kinetic models of economy}\label{graz2} In
this section we show how to extend the asymptotic analysis of the
previous section to the case in which the kinetic model describes
the time evolution of a density  $f(v,t)$, which now denotes the
distribution of wealth $v \in \real_+$ among economic agents at
time $t \ge 0$. The collision \fer{coll} represents now a trade
between individuals. For a deep insight into the matter, we
address the interested reader to \cite{YD, GGPS, IKR, PGS, BP},
and to the references therein. With the convention $f(v,t) = 0$ if
$v < 0$, the kinetic model reads \cite{CPT, P}.
 \be
 \frac{\partial f(v)}{\partial t} =
\int_{\real_+}\left( \frac 1J f(v_*)
 f(w_*) - f(v) f(w)\right) d
 w \label{eq:bolt1}
  \ee
where $(v^*,w^*)\in \real_+$ are the post-trade wealths  generated
by the couple $(v,w)$ after the interaction, along the rule
\fer{coll}. As before, the jacobian $J = p^2-q^2$. Since the
$v$-variable takes values in $\real_+$, the collision rules
\fer{coll} lead to a remarkable difference with respect to the
case treated in the previous section. The pair $(v_*,w_*)$ of
pre-collision variables that generate the pair $(v,w)$ is given by
\[
v_* = \frac{pv -qw}J , \quad w_* = \frac{pw -qv}J.
\]
While in the former case this pair is always admissible ( $v_*,w_*
\in \real$), in the latter we have to discard all pairs of
pre-collision variables for which $v_* <0$ or $w_* <0$. This shows
that, for any given $v \in \real_+$,   the product $f(v_*)
 f(w_*)$ in \fer{eq:bolt1} is different from zero only on the set $\bigB = \{(q/p)v < w <
 (p/q)v\}$. This implies in other words that, if we fix the wealth $v \in \real_+$
 as outcome of a single trade, the other outcome $w$ can only lie
 on the subset $\bigB$.

A great simplification is obtained writing equation \fer{eq:bolt1}
 in weak form, where the presence of pre-collision wealths is
 avoided,
 \be\label{eco}
\frac d{dt}\int_{\real_+} \phi(v)f(v,t)\,dv   =
 \int_{\real_+^2} f(v,t)f(w,t) ( \phi(v^*)-\phi(v)) dv d w .
 \ee

\begin{remark}
The role of the energy is now played by the mean $m(t) = \int
vf(v,t)\, dv$. Note however that one can think to equation
(\ref{eq:bolt1}) as the analogous of the isotropic form of a
hard-sphere Boltzmann equation for a density function $f(v',t)$,
$v'\in\real$ written with respect to energy variable $v=(v')^2/2$.
In this sense it is again the non conservation of the energy that
will originate the power law tails.
\end{remark}

To look for self--similarity we scale our solution according to
 \be
  g(v,t) = { m(t)}f\left( { m(t)}v, t \right),
  \ee
 which implies that $\int v g(v,t)= 1$ for all $t \ge 0$. Moreover $g= g(v,t)$ satisfies
 the equation
 \[\nonumber
\frac d{dt}\int_{\real_+} \phi(v)g(v,t)\,dv  -\left(p +q
-1\right)\int_{\real_+}\phi(v)\frac{\partial}{\partial v}\left(vg\right)\, dv =
 \]
 \be\label{weak boltz11}
 \int_{\real_+^2} g(v)
g(w) ( \phi(v^*)-\phi(v)) dv d w .
 \ee
Performing the same computations of the previous section, and {\it mutatis mutandis} we conclude
that $g(v,t)$ satisfies
 \[
\frac d{dt}\int_{\real_+} \phi(v)g(v,t)\,dv + q\int_{\real}\phi'(v)(v -1)g(v)\, dv =
 \]
 \be
\frac 12\int_{\real} g(v)\left((p-1)^2v^2 +q^2w^2 +2(p-1)qvw\right)\phi''(v) dv + R(p,q).
 \ee
 The form of the remainder $R(p,q)$ is analogous to that of \fer{resto}. It is clear that the
 correct scaling for small values of the parameter $q$ is now
 \be\label{resc2}
 \tau = q t, \quad h(v,\tau) = g(v,t),
 \ee
 which implies that $h(v, \tau)$ satisfies the equation
 \[
\frac d{d\tau}\int_{\real_+} \phi(v)h(v,\tau)\,dv + \int_{\real}\phi'(v)(v -1)h(v)\, dv =
 \]
 \be\label{quasiFF}
\frac 12\int_{\real} h(v)\frac{(p-1)^2}q v^2\phi''(v) dv     + R_1(p,q),
 \ee
 where the remainder $R_1$ is given by
 \[
 R_1(p,q) = \frac 12\int_{\real_+}\left( qw^2 +2(p-1)vw\right)\phi''(v) dv + \frac
1qR(p,q).
 \]
Let us consider a parameter $p=p(q)$ such that, for a given
constant $\la>0$
  \be\label{con2}
 \lim_{q\to 0} \frac{(p(q)-1)^2}q = \lambda.
  \ee
 Then,  equation \fer{quasiFF} is well--approximated by the equation (in weak form)
 \be\label{FF}
 \frac d{d\tau}\int_{\real} \phi(v)h(v,\tau)\,dv +\int_{\real}\phi'(v)(v-1) h(v)\, dv =
\frac\lambda{ 2}\int_{\real} h(v)v^2\phi''(v) dv .
 \ee
 Equation \fer{FF} is nothing but the weak form of the Fokker-Planck equation
 \be\label{FP2b}
 \frac{\partial h}{\partial \tau} = \frac \lambda{2}\frac{\partial^2 }
 {\partial v^2}\left( v^2 h\right) + \frac{\partial }{\partial v}\left(v
 h\right),
 \ee
 which admits a unique stationary state of unit mass, given by the
 $\Gamma$-distribution \cite{BM, CPT}
 \be\label{equi}
M_\lambda(v)=\frac{(\mu-1)^\mu}{\Gamma(\mu)}\frac{\exp\left(-\frac{\mu-1}{v}\right)}{v^{1+\mu}}
 \ee
  where
  $$ \mu = 1 + \frac{2}{\lambda} >1.
$$
This stationary distribution exhibits a Pareto power law tail for
large $v$'s.

Note that this equation is essentially the same Fokker-Planck equation derived from a
Lotka-Volterra interaction in \cite{BM, So, BMRS}.

\begin{remark}
The formal analysis shows that the Fokker--Planck equation
\fer{FF} follows from the kinetic model independently of the sign
of the quantity $p+q-1$, which can produce exponential growth of
wealth (when positive), or exponential dissipation of wealth (when
negative). Hence, Pareto tails are produced in both situations, as
soon as the compatibility condition \fer{con2} holds. As discussed
in Remark \ref{poss}, condition \fer{con2} is always admissible if
$p+q-1 <0$, while  one has to require $p
>1$ if $p+q-1 >0$. This is quite remarkable
since it shows that this uneven distribution of money which
characterizes most western economies may not only be produced as
the effect of a growing economy but also under critical economical
circumstances.
\end{remark}

\begin{remark}
The model studied in \cite{Sl} corresponds to the choice $p = 1-q +\e$, with $\e >0$. This
interaction implies exponential growth of wealth, and convergence of the solution to the
Fokker-Planck equation if $\e = \e(q)$ satisfies
 \[
 \lim_{q \to 0} \frac{\e^2(q)}q = \lambda.
 \]
 Since the same limit equation is derived within the choice $p = 1-q -\e$, we are free to choose
 $\e$ negative. The particular choice
  \[
  \e = -2\sqrt q +2q,
  \]
 which implies $\mu = 3/2$ and thus $\lambda=4$, leads to the stationary state \cite{Sl}
  \be
  \label{eq:std3}
 M_{4}(v)=\frac1{\sqrt{2\pi}}\frac{\exp\left(-\frac{1}{2v}\right)}{v^{5/2}},
 \ee
 which solves the kinetic equation \fer{weak boltz11} for all values of the scaling parameter
 $q<1/4$.
\end{remark}

\section{The Fourier transform of the kinetic equations} \label{appl}

The formal results of Sections \ref{graz1} and \ref{graz2} suggest
that, at least in the limit $p \to 1$ and $q \to 0$, the
large-time behavior of the solution to the kinetic model
\fer{boltz1} is characterized by the presence of overpopulated
tails. In what follows, we will justify rigorously this behavior,
at least for a certain domain of the mixing parameters $p$ and
$q$. We start our analysis with a detailed study of the Boltzmann
model \fer{eq:boltz}.

The initial value problem for this model can be easily studied
using its weak form \fer{weak boltz}.
 Let $\bigM_0$ the space of all probability measures in $\real_+$ and by
  \bq\label{misure} \bigM_{\a} =\left\{
\mu \in\bigM_0: \int_{\real} |v|^{\a}\mu(dv) < +\infty, \a\ge
0\right\},
 \eq
 the space of all
Borel probability measures of finite momentum of order $\a$,
equipped with the topology of the weak convergence of the
measures.

 By a weak solution of the initial value problem for equation
\rf{eq:boltz}, corresponding to the initial probability density
$f_0(w) \in \bigM_{\a}, \a >2$ we shall mean any probability
density $f \in C^1(\real, \bigM_{\a})$ satisfying the weak form of
the equation
 \be\label{weakk}
\frac d{dt}\int_{\real} \phi(v)f(v,t)\,dv  =
 \int_{\real^2} f(v)
f(w) ( \phi(v^*)-\phi(v)) dv d w ,
 \ee
 for $t>0$ and all smooth functions $\phi$, and such
that for all $\phi$
 \bq\label{ic} \lim_{t\to 0} \int_{\real} \phi(v)f(v,t)\, dv = \int_{\real} \phi(v)f_0(v)\, dv.
  \eq
In the rest of this section, we shall study the weak form of
equation \rf{eq:boltz}, with the normalization conditions
\fer{norm1}. It is equivalent to use the Fourier transform of the
equation~\cite{Bob}:
\begin{equation} \label{eq}
\frac{\partial \widehat f(\xi, t)}{\partial t}  =
\widehat{Q}\left( \ff,\ff \right)(\xi,t) ,
\end{equation}
where $\ff(\xi,t)$ is the Fourier transform of $f(x,t)$,
\[ \widehat{f}(\xi,t) = \int_{{\real}} e^{-i \xi v}\, f(v,t)\, dv, \]
and
\begin{equation}\label{trascoll}
  \widehat{Q}\left( \ff,\ff \right)(\xi) =
 \widehat{f}(p\xi)
\widehat{f}(q\xi) - \widehat{f}(\xi) \widehat{f}(0).
\end{equation}
The initial conditions \fer{norm1} turn into
\[ \widehat{f}(0)=1,  \widehat{f}^\prime(0) =0, \widehat{f}^{\prime\prime}(0) = -1, \]
$\widehat{f} \in C^2(\real)$. Hence equation \fer{eq} can be
rewritten as
\begin{equation}\label{fkac}
\frac{\partial \widehat f(\xi, t)}{\partial t}  + \ff(\xi,t) =
\widehat{f}(p\xi) \widehat{f}(q\xi) .
\end{equation}
Equation \fer{fkac} is a special case of equation (4.8) considered
by Bobylev and Cercignani in \cite{BC4}. Consequently, most of
their conclusions applies to the present situation as well.  The
main difference here is that the mixing parameters $p$ and $q$ are
allowed to assume values bigger than $1$.

 We introduce a metric on $\bigM_{p}$ by
\begin{equation} \label{ds}
d_s (f,g) = \sup_{\xi \in {\real}} \frac{|\widehat{f}(\xi) -
\widehat{g}(\xi)|}{|\xi|^s}
\end{equation}

Let us write $s = m+\alpha$, where $m$ is an integer and $0 \leq
\alpha < 1$. In order that $d_s(F,G)$ be finite, it suffices that
$F$ and $G$ have the same moments up to order $m$.

The norm \fer{ds} has been introduced in~\cite{GTW} to investigate
the trend to equilibrium of the solutions to the Boltzmann
equation for Maxwell molecules. There, the case $s = 2 + \alpha$,
$\alpha > 0$, was considered. Further applications of $d_s$ can be
found  in~\cite{CGT, CCG1, TV1, GJT}.


\subsection{Uniqueness and asymptotic behavior}\label{uni}

We will now study in details the asymptotic behavior of the scaled
function $g(v,t)$. As briefly discussed before, a related analysis
has been performed in the framework of the study of self--similar
profiles for the Boltzmann equation for Maxwell molecules in
\cite{BC4, BCT}. Likewise, the role of the Fourier distance in the
asymptotic study of nonconservative kinetic equations has been
evidenced  in \cite{PT}. Consequently, part of the results
presented here fall into the results of \cite{BC4, PT}, and could
be skipped. Nevertheless, for the sake of completeness, we will
discuss the point in an exhaustive way.

 The
existence of a solution to equation \fer{eq:boltz} can be seen
easily using the same methods available for the elastic Kac model.
In particular, a solution can be expressed as a Wild sum
\cite{Bob, CGT}. In order to prove uniqueness, we use the method
first introduced in~\cite{GTW}.  Let $f_1$ and $f_2$ be two
solutions of the Boltzmann equation \fer{eq:boltz},
 corresponding to initial values $f_{1,0}$ and $f_{2,0}$ satisfying
 conditions \fer{norm1},
and $\widehat{f_1}$, $\widehat{f_2}$ their Fourier transforms.
Given any positive constant $s$, with $2 \le s \le 3$, let us
suppose in addition that $d_s(f_{1,0},f_{2,0})$ is bounded. Then,
it holds
 \begin{equation} \label{ineqd}
\frac{\partial}{\partial t} \frac{ \left ( \widehat{f_1} -
\widehat{f_2} \right )}{|\xi|^s} + \frac{ \widehat{f_1}(\xi)-
\widehat{f_2}(\xi)}{|\xi|^s} =
 \frac{\widehat{f_1}(p \xi)\widehat{f_1}(q \xi) -
\widehat{f_2}(p \xi)\widehat{f_2}(q \xi)}{|\xi|^s}.
 \end{equation}
Now, since $|\widehat{f_1}(\xi)| \le 1$ ($|\widehat{f_2}(
 \xi)| \le 1$), we obtain
 \be\nonumber
 \left | \frac{\widehat{f_1}(p \xi)\widehat{f_1}(q \xi) - \widehat{f_2}(p
 \xi)\widehat{f_2}(q \xi)}{|\xi|^s} \right | \leq |\widehat{f_1}(p
 \xi)| \left | \frac{\widehat{f_1}(q \xi)-\widehat{f_2}(q
 \xi)}{|q \xi|^s} \right |q^s
 +
 \ee
 \be
 + |\widehat{f_2}(q \xi)| \left | \frac{\widehat{f_1}(p
 \xi)-\widehat{f_2}(p \xi)}{|p \xi|^s} \right |p^s
 \leq \sup \left | \frac{\widehat{f_1} -
\widehat{f_2}}{|\xi|^s} \right|(p^s+q^s).
 \ee
We set
\[ h(t,\xi) = \frac{\widehat{f_1}(\xi) - \widehat{f_2}(\xi)}{|\xi|^{s}}. \]

The preceding computation shows that
\begin{equation} \label{estim}
\left| \frac{\partial h}{\partial t} + h \right| \leq  (p^s +q^s)
\|h\|_\infty .
\end{equation}
Gronwall's lemma proves at once that
 \[
\|h(t)\|_{\infty} \le \exp\left\{ (p^s + q^s -1)t\right\}
\|h_0\|_{\infty} .
 \]
We have

\begin{theorem} \label{contr}
Let $f_1(t)$ and $f_2(t)$ be two solutions of the Boltzmann
equation~\fer{eq:boltz}, corresponding to  initial values
$f_{1,0}$ and $f_{2,0}$ satisfying conditions \fer{norm1}. Then,
if for some  $2 \le s \le 3$,  $d_s(f_{1,0},f_{2,0})$ is
 bounded, for all times $t \geq 0$,
 \be\label{dec3}
  d_s(f_1(t), f_2(t)) \leq \exp\left\{ (p^s + q^s -1)t\right\}d_s(f_{1,0},f_{2,0}).
  \ee
In particular, let $f_0$ be a nonnegative density satisfying
conditions \fer{norm1}. Then, there exists a unique weak solution
$f(t)$ of the Boltzmann equation, such that $f(0) = f_0$. In case
$p^s + q^s -1 <0$ the distance $d_s$ is contracting exponentially
in time.
\end{theorem}

\begin{figure}[ht]
\begin{center}
\includegraphics[scale=.5]{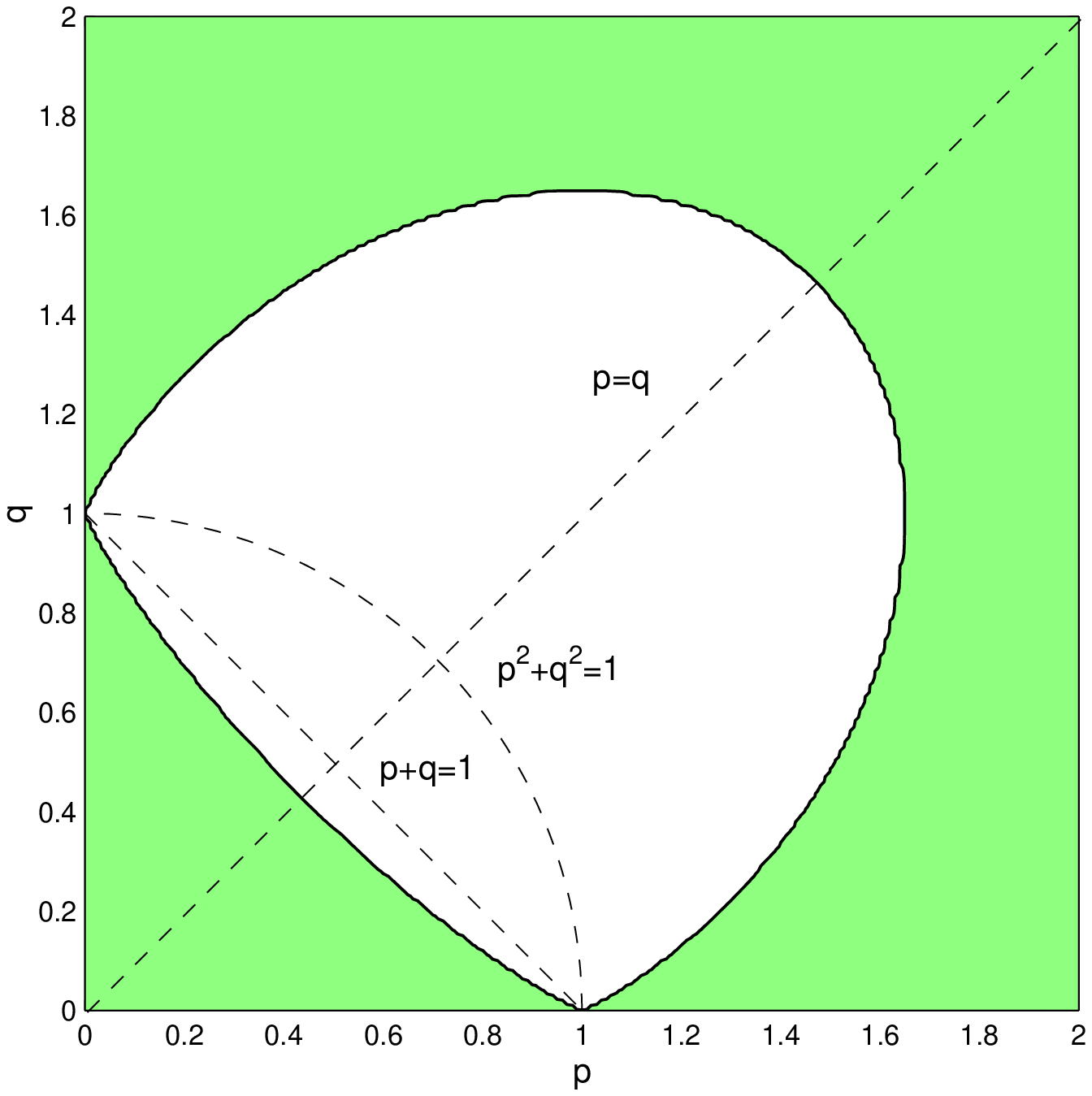}
\end{center}
\caption{The white domain represents the region where the minimum
of the function $\SS$ is negative for $p,q \in [0,2]$.}
\label{fig:0}
\end{figure}

Let us remark that, given a constant $a>0$,
 \begin{equation} \label{ds1}
 \sup_{\xi \in {\real}} \frac{|\widehat{f_1}(a\xi) -
\widehat{f_2}(a\xi)|}{|\xi|^s} =  a^s \sup_{\xi \in {\real}}
\frac{|\widehat{f_1}(a\xi) - \widehat{f_2}(a\xi)|}{|a\xi|^s}= a^s
d_s (f_1,f_2).
 \end{equation}
Hence, if $g(t)$ represents the solution $f(t)$ scaled by its
energy like in \fer{scala},
 \[
 \widehat{g}(\xi) = \widehat{f}\left(\frac\xi{\sqrt{
 E(t)}}\right),
 \]
and from \fer{ds1} we obtain the bound
 \begin{equation} \label{ds11}
 d_s (g_1(t),g_2(t))= \sup_{\xi \in {\real}} \frac{|\widehat{g_1}(\xi,t) -
\widehat{g_2}(\xi,t)|}{|\xi|^s} =  \left(\frac 1{\sqrt{
 E(t)}}\right)^s d_s (f_1(t),f_2(t)).
 \end{equation}
 Using \fer{dec3}, we finally conclude that, if $g_1(t)$ and
 $g_2(t)$ are two solutions of the scaled Boltzmann
equation~\fer{boltz1}, corresponding to  initial values $f_{1,0}$
and $f_{2,0}$ satisfying conditions \fer{norm1},
 then, if $2 \le s \le 3$,
for all times $t \geq 0$,
 \be\label{dec1}
  d_s(g_1(t), g_2(t)) \leq \exp\left\{\left[ (p^s + q^s -1)
   -\frac s2(p^2+q^2 -1)\right]t\right\}d_s(f_{1,0},f_{2,0}).
  \ee
Let us define, for $\d \ge 0$,
 \be\label{key}
 \SS (\delta) = p^{2+\delta} + q^{2+\delta} -1 -\frac{2+\d}2\left(p^2+q^2
 -1\right).
 \ee
Then, the sign of $\SS$ determines the asymptotic behavior of the
distance $d_s(g_1(t), g_2(t))$. In particular, if there exists an
interval $0 < \d < \bar\d$ in which $\SS (\delta) <0$, we can
conclude that $d_{2+\d}(g_1(t), g_2(t))$ converges exponentially
to zero. Note that, by construction, $\SS (0)= 0$, and thus
$\min_\delta\{\SS\}\leq 0$. The function \fer{key} was first
considered by Bobylev and Cercignani in \cite{BC4}. The sign of
$\SS$, however was studied  mainly for $p = 1-q$, namely the case
of the dissipative Boltzmann equation. In Figure \ref{fig:0} a
numerical evaluation of the region where the minimum of the
function $\SS$ is negative for $p,q \in [0,2]$ is reported.

\begin{remark}
The behavior of $\SS (\delta)$ when $p^2+q^2 =1$ is clear. In this
case in fact, both $p$ and $q$ are less than $1$, which implies
 \[
\SS (\delta) = p^{2+\delta} + q^{2+\delta} -1 <0,
 \]
 for all $\d >0$. We can draw the same conclusion when $p^2+q^2
 >1$, while both $p<1$ and $q<1$.

 Consider now the case $p^2+q^2 > 1$, with $p>1$. In this case,
  while $-\frac{2+\d}2\left(p^2+q^2 -1\right)$ decreases linearly,
$p^{2+\d}$ increases exponentially, and the sign of $\SS(\d)$
becomes positive  for large values of $\d$.

 If finally  $p^2+q^2 < 1$, the sign of $\SS(\d)$ for large values
 of $\d$ is positive, since, while $p^{2+\delta} + q^{2+\delta} -1
 \ge -1$,
 \[
 -\frac{2+\d}2\left(p^2+q^2 -1\right) \ge 1
  \]
  for
  \[
  \d \ge \frac{2(p^2+q^2)}{1 -(p^2+q^2)}.
  \]

\end{remark}

The previous remark indicates that in the general case one can at
best  hope that $\SS(\d)$ is negative  in an interval $(0,
\bar\d)$. To show that this is really the case, one has to
investigate carefully the behavior of $\SS(\d)$ in a neighborhood
of zero. Since the function $\SS(\d)$ is convex for $\d \ge 0$,
 \[
 \frac{d^2\SS(\d)}{d\d^2} = p^{2+\delta}(\log p)^2 +
 q^{2+\delta}(\log q)^2 >0 ,
 \]
 and $\SS(0)= 0$, in all cases where $\SS(\d)$ is positive for
 large values of $\d$, a sufficient condition for $\SS(\d)$ be
 negative in some interval $0 < \d < \bar\d$ is that
 \[
 \left. \frac{d\SS(\d)}{d\d}\right|_{\d=0} = p^{2}\log p +
 q^{2}\log q - \frac{1}2\left(p^2+q^2
 -1\right) <0 .
 \]

 Let us discuss before the case $p^2+q^2 < 1$. Given $\la >0$, we introduce a dependence
 between $p$ and $q$ by setting $p = 1 -\la q$. Since $p >q$, this
 relationship is possible only if
  $q< 1/(1+\la)$. Moreover $p^2+q^2 < 1$ requires $q <
 (2\la)/(1+\la^2)$. Using this, it is immediate to show that there is an interval $0\le
 q \le \bar q$ in which $\SS'(0) <0$. We have
 \[
 \left. \frac{d\SS(\d)}{d\d}\right|_{\d=0} = G(q) = (1-\la q)^{2}\log (1-\la q) +
 q^{2}\log q - \frac{1}2\left((1-\la q)^2+q^2
 -1\right).
 \]
 Clearly, $G(0) = 0$. Moreover
  \[
  G'(q) = 2q\log q -2\la(1-\la q)\log(1-\la q),
  \]
  and
  \[
  G''(q) = 2(1 +\log q) +2\la^2\left( 1 + \log(1-\la q)\right).
  \]
  Now $G''(q)<0$ in some interval $(0, q_1)$, which implies that
  $G'(q)$ is decreasing in the same interval. But, since
  $G'(0)=0$, $G'(q) <0 $ in the interval $(0, \bar q)$, where $\bar
  q$ solves
  \[
   2\bar q\log \bar q -2\la(1-\la\bar q)\log(1-\la \bar q)= 0
  \]
 Consequently, $G(q) <0$ at least in the same interval.

 Let us now treat the case $p^2+q^2 > 1$, with $p>1$. Let us set $p = 1 +\la
 q$.
  We have
 \[
 \left. \frac{d\SS(\d)}{d\d}\right|_{\d=0} = G(q) = (1+\la q)^{2}\log (1+\la q) +
 q^{2}\log q - \frac{1}2\left((1+\la q)^2+q^2
 -1\right).
 \]
 In this case
  \[
  G'(q) = 2q\log q +2\la(1+\la q)\log(1+\la q),
  \]
  and
  \[
  G''(q) = 2(1 +\log q) +2\la^2\left( 1 + \log(1+\la q)\right).
  \]
  As before, $G''(q)<0$ in some interval $(0, q_2)$, which implies that
  $G'(q)$ is decreasing in the same interval. But, since
  $G'(0)=0$, $G'(q) <0 $ in the interval $(0, \bar q)$, where $\bar
  q$ now solves
  \[
   2\bar q\log \bar q + 2\la(1+\la\bar q)\log(1+\la \bar q)= 0
  \]
 Consequently, $G(q) <0$ at least in the same interval.

 We proved

 \begin{lemma}\label{SS}
 Let $\SS(\d), \d \ge 0$ be the function defined by \fer{key}.
 Given a constant $\la >0$, if $p^2+q^2 < 1$, let us define $p=1- \la q$ . Then,
 provided $q < \min\left\{1/(1+\la),(2\la)/(1+\la^2)\right\}$ there exists an interval $I_-= (0 , \bar \d_-(q))$
 such that $\SS(\d)<0$ for $\d \in I_-$. If $p^2+q^2 > 1$, and $p=1+ \la
 q$ there exists an interval $I_+= (0 , \bar \d_+(q))$
 such that $\SS(\d)<0$ for $\d \in I_+$. In the remaining cases, namely when $p^2+q^2 =
 1$ or $p^2+q^2 > 1$ but $p<1$, $\SS(\d)<0$ for all $\d>0$.
 \end{lemma}

Lemma \ref{SS} has important consequences both in the behavior of
the solution to the Boltzmann equation \fer{boltz1}, and in the
limit procedure introduced in Sections \ref{graz1} and
\ref{graz2}. The main consequence of the lemma is contained into
the following.

\begin{theorem} \label{de}
Let $g_1(t)$ and $g_2(t)$ be two solutions of the Boltzmann
equation~\fer{boltz1}, corresponding to  initial values $f_{1,0}$
and $f_{2,0}$ satisfying conditions \fer{norm1}. Then, there
exists a constant $\bar\d >0$ such that, if $2 < s < 2+ \bar\d$,
for all times $t \geq 0$,
 \be\label{decc}
  d_s(g_1(t), g_2(t)) \leq \exp\left\{ -C_st\right\}d_s(f_{1,0},f_{2,0}).
  \ee
 The constant $C_s = -\SS(s-2)$ is strictly positive, and
 the distance $d_s$ is contracting exponentially
in time.
\end{theorem}


\subsection{Convergence to self--similarity}\label{convv}

By means of the estimates of Section \ref{uni}, we will now
discuss the evolution of moments for the solution to equation
\fer{boltz1}. By construction, the second moment of $g(v,t)$ is
constant in time, and equal to $1$ thanks to the normalization
conditions \fer{norm1}. We can  use the computations leading to
the Fokker-Planck equation \fer{FP2}, choosing $\phi(v) =
|v|^{2+\d}$, where for the moment the positive constant $\d \le
1$. Suppose  that the initial density $g_0(v) = f_0(v)$ is such
that
\begin{equation}\label{moment}
\int_{\real } |v|^{2+\d}g_0(v)\, dv = m_\delta < \infty.
\end{equation}
Then, since the contribution due to the term $
\frac\partial{\partial v}{\left(vg(v)\right)}$ can be evaluated
integrating by parts,
 \[
\int_{\real} |v|^{2+\delta} \frac\partial{\partial v}{\left(
vg(v)\right)}\, dv = -(2+\delta)\int_{\real}|v|^{2+\delta} g(v,t)
\, dv ,
\]
we obtain
\begin{equation}\nonumber
\displaystyle \frac d{dt}\int_{\real}|v|^{2+\delta}g(v,t)\, dv
+(2+\delta)
 \g  \int_{\real}|v|^{2+\delta}g(v,t) \, dv =
\end{equation}
\begin{equation}\label{inv1}
\int_{\real^2}\,dv\,dw \left(
|pv+qw|^{2+\delta}-|v|^{2+\delta}\right)g(v) g(w) \, .
\end{equation}
Let us recover a suitable upper bound for the last integral in
\fer{inv1}. Given any two constants $a, b$, and $0< \d \le 1$ the
following inequality holds
 \begin{equation}\label{in1}
\left( |a| +|b| \right)^\d \le |a|^\d + |b|^\d .
 \end{equation}
Hence, choosing $a =p|v|$ and $b=q|w|$,
 \[
|pv+qw|^{2+\d} \le \left( pv+qw \right)^2\left(p^\d|v|^\d +
q^\d|w|^\d\right).
 \]
Substituting into the right-hand side of \fer{inv1}, recalling
that the mean value of $g$ is equal to zero, and the second moment
of $g$ equal to one, gives
 \begin{equation}\nonumber
\int_{\real^2} |pv+qw|^{2+\d}g(v) g(w) \, \,dv\,dw \le
 \end{equation}
  \begin{equation}\nonumber
\int_{\real^2} (pv+qw)^{2}\left(p^\d|v|^d + q^\d|w|^d\right)g(v)
g(w) \, \,dv\,dw =
 \end{equation}
   \begin{equation}\nonumber
\left( p^{2+\d}+ q^{2+\d}\right)\int_{\real} |v|^{2+\delta}g(v)
 \,dv\, + \left(p^{2}q^\d+ q^{2}p^\d\right)\int_{\real} |v|^{\d}\,
 dv .
 \end{equation}

Grouping all these inequalities, and recalling the expression of
$\SS(\delta)$ given by \fer{key} we obtain the differential
inequality
\begin{equation}\label{dd}
\displaystyle \frac d{dt}\int_{\real}|v|^{2+\delta}g(v,t)\, dv \le
\SS(\delta)\int_{\real}|v|^{2+\delta}g(v,t) \, dv + B_{p,\delta},
\end{equation}
where, by H\"older inequality
 \be\label{Bpl}
 B_{p,\delta} \le p^{2}q^\d+ q^{2}p^\d.
\ee
 By Lemma \ref{SS}, for any $\delta <\bar\delta$, $\SS(\delta) <0$. In this
case, inequality \fer{dd} gives an upper bound for the moment,
that reads
\begin{equation}\label{mom}
\int_{\real}|v|^{2+\delta}g(v,t)\, dv \le m_\delta +
\frac{B_{p,\delta}}{|\SS(\delta)|} < \infty.
\end{equation}
 In the case $\bar\delta >3$ we can easily iterate our procedure to obtain that any moment of order
$2+ \delta$, with $\delta < \bar\delta$ which is  bounded
initially, remains bounded at any subsequent time. The only
difference now is that the explicit expression of the bound is
more and more involved.

If $\delta < \bar\delta$, we can immediately draw conclusions on
the large--time convergence of class of probability densities
$\left\{ g(v, t)\right\}_{t \ge 0}$,  By virtue of Prokhorov
theorem (cfr. \cite{LR}) the existence of a uniform bound on
moments implies that this class is tight, so that any sequence
$\left\{ g(v, t_n)\right\}_{n \ge 0}$ contains an infinite
subsequence which converges weakly to some probability measure
$g_{\infty}$. Thanks to our bound on moments, provided $\d < \bar
\d$, $g_{\infty}$ possesses moments of order $2+\d$, for $0 <
\delta < \bar\delta$.

It is now immediate to show that this limit is unique. To this
aim, let us consider two initial densities $f_{0,1}(v)$ and
$f_{0,2}(v)$ such that, for some $0 <\delta <\bar\delta$,
\[
\int_{\cal R}|v|^{2+\d}f_{0,1}(v)\,dv < +\infty, \quad \int_{\cal
R}|v|^{2+\d}f_{0,2}(v)\,dv < +\infty.
\]
Then, by Theorem \ref{de}, the distance $d_s(f_1(t),f_2(t))$
between the solutions converges exponentially to zero with respect
to time, as soon as $2 < s < 2+ \bar\d$. Let now $f_0(v)$ possess
finite moments of order $2+ \d$, with $0 <\delta <\bar\delta$.
Thanks to our previous computations on moments, for any fixed time
$T>0$, the corresponding solution $f(v,T)$ has finite moments of
order $2+ \d$. Choosing $f_{0,1}(v)= f_0(v)$, and $f_{0,2}(v)=
f(v,T)$ shows that $d_s(f(t),f(t+T))$ converges exponentially to
zero in time. It turns out that the $d_{s}$-distance between
subsequences  converges to zero as soon as $\SS(s-2) <0$.

We can now show that the limit function $g_{\infty}(v)$ is a
stationary solution to \fer{boltz1}. We know that if condition
\fer{moment} holds,  both the solution $g(v,t)$ to equation
\fer{boltz1} and $g_{\infty}(v)$ have moments of order $2+
\delta$, with $0 < \delta < \bar\delta$ uniformly bounded. Hence,
for any $t \ge 0$, proceeding as in the proof of Theorem
\ref{contr}, we obtain
\begin{equation}
d_{s}\left( Q(g(t),g(t)), Q(g_{\infty},g_{\infty})\right) \le
\left( p^s+q^s +1 \right) d_{s}\left( g(t), g_{\infty}\right).
\end{equation}
This implies the weak* convergence of $Q(g(t),g(t))$ towards
$Q(g_{\infty},g_{\infty})$. In particular, due to the equivalence
among different metrics which metricize the weak* convergence of
measures \cite{GTW, TV1}, if $C^1_0({\R})$ denotes the set of
compactly supported continuously differentiable functions, endowed
with its natural norm $\| \cdot \|_1$, for all $\phi \in
C^1_0({\R})$,
\begin{equation}
\int_{\real}\phi(v)Q(g(t),g(t))(v) \, dv \to
\int_{\real}\phi(v)Q(g_{\infty},g_{\infty})(v) \, dv.
\end{equation}
On the other hand, for all $\phi \in C^1_0({\R})$, integration by
parts gives
\begin{equation}
\int_{\real}\phi(v)\frac{\partial}{\partial v}\left(vg(v,t)\right)
\, dv = - \int_{\real}v\phi^{\, \prime} (v)g(v, t)\, dv.
\end{equation}
Since $|v\phi^{\, \prime} (v)| \le |v|\| \phi^{\, \prime} \|_1$,
and the second moment of $g(v,t)$ is equal to unity, the
convergence of $d_{s}\left( g(t), g_{\infty}\right)$ to zero
implies
\begin{equation}
\int_{\real}v\phi^{\, \prime} (v)g(v, t) \, dv \to
\int_{\real}v\phi^{\, \prime} (v)g_{\infty}(v) \, dv.
\end{equation}
Finally, for all $\phi \in C^1_0({\R})$ it holds
\begin{equation}
\int_{\real}\phi(v)\left\{ \frac{\partial}{\partial
v}\left(vg_{\infty}(v)\right) -
Q(g_{\infty},g_{\infty})(v)\right\} \, dv = 0 .
\end{equation}
This shows that $g_{\infty}$ is the unique stationary solution to
\fer{boltz1}. We have
 \

\begin{theorem} \label{pnot1}
Let $\d >0$ be such that $\SS(\d)<0$, and let $g_{\infty}(v)$ be
the unique stationary solution to equation \fer{boltz1}. Let
$g(v,t)$ be the weak solution of the Boltzmann
equation~\fer{boltz1}, corresponding to the initial density $f_0$
satisfying
\[
\int |v|^{2+\delta} \, f_0(v)\, dv<\infty .
\]
 Then, $g(v,t)$ satisfies
\[
\int |v|^{2+\delta} \, g(v, t)\, dv \le c_{\delta} <\infty .
\]
If $0< \delta\le 1$ the constant $c_{\delta}$ is given by
\fer{mom}. Moreover,  $g(v,t)$ converges exponentially fast in
Fourier metric towards $g_{\infty}(v)$, and the following bound
holds
\begin{equation}\label{conv1}
d_{2+\delta}(g(t),g_{\infty} ) \le d_{2+\delta}(f_0,g_{\infty}
)\exp\left\{-|\SS(\delta)|t\right\}
\end{equation}
where $\SS(\delta)$ is given by \fer{key}.
\end{theorem}

Depending of the values of the mixing parameters $p$ and $q$, the
stationary solution $g_{\infty}$ can have overpopulated tails. We
can easily check the presence of overpopulated tails by looking at
the singular part of the Fourier transform \cite{EB1}. Since the
Fourier transform of $g_{\infty}$ satisfies the equation
\begin{equation}\label{foug}
-\g\xi \frac{\partial{\widehat{g}}}{\partial \xi}+\widehat
g(\xi)=\widehat g(p\xi)\widehat g(q\xi),
\end{equation}
we set
 \be\label{serie}
\widehat{g}(\xi) = 1-|\xi|^2 + A|\xi|^{2+\delta} + \dots
 \ee
 which takes into account the fact that  $g_{\infty}$ satisfies
 conditions \fer{norm1}. The leading small $\xi$-behavior of the
 singular component will reflect an algebraic tail of the velocity
 distribution. Substitution of expression \fer{serie} into
 \fer{foug} shows that the coefficient of the power
 $|\xi|^{2+\delta}$ is $A\SS(\d)$. Thus, the term
 $A|\xi|^{2+\delta}$ can appear in the expansion of
 $\widehat{g}(\xi)$ as soon as $\d$ is such that $\SS(\d) = 0, \d >0$.
In other words, tails in the stationary distributions are present
in all cases in which there exists a $\d = \bar \d>0$ such that
$\SS(\bar\d)=0$. Now the answer is contained into Lemma \ref{SS}.


\subsection{The grazing collision asymptotics}\label{grazz}

The results of the previous section are at the basis of the
rigorous derivation of the Fokker-Planck asymptotics formally
derived in Section \ref{graz1} and \ref{graz2}. Suppose  that the
initial density $g_0(v) = f_0(v)$ satisfies condition
\fer{moment}. Using a Taylor expansion, we obtain
 \[
 |pv+qw|^{2+\delta}-|v|^{2+\delta} =
 \]
\be
 (2+\d)|v|^\d v((p-1)v+qw) +
 \frac 12 (1+\d)|\tilde v|^\d((p-1)v+qw)^2,
 \ee
where, for some $0\le \theta \le 1 $,
 \[
 \tilde v = \theta (pv+qw) +(1-\theta)v.
 \]
Using this into  equality \fer{inv1}, one has
\begin{equation}\nonumber
\displaystyle \frac d{dt}\int_{\real}|v|^{2+\delta}g(v,t)\, dv
+(2+\delta)
 \g  \int_{\real}|v|^{2+\delta}g(v,t) \, dv =
\end{equation}
 \[
(2+\d)\int_{\real^2} \left(|v|^\d v((p-1)v+qw) \right)g(v)
g(w)\,dv\,dw
 \]
\begin{equation}\label{inv3}
 + \frac 12 (2+\d)(1+\d) \int_{\real^2}|\tilde
v|^\d((p-1)v+qw)^2 g(v) g(w)\,dv\,dw.
\end{equation}
Since the momentum of $g$ is equal to zero, we can rewrite
\fer{inv3} as
\begin{equation}\nonumber
\displaystyle \frac d{dt}\int_{\real}|v|^{2+\delta}g(v,t)\, dv
+\frac{2+\delta}2\left[ (p-1)^2 + q^2 \right]
  \int_{\real}|v|^{2+\delta}g(v,t) \, dv \le
\end{equation}
\begin{equation}\label{inv4}
 + \frac 12 (2+\d)(1+\d) \int_{\real^2}|\tilde
v|^\d((p-1)v+qw)^2 g(v) g(w)\,dv\,dw.
\end{equation}
Assuming $0 <\d <1$,
\[
 |\tilde v| \le (1+p)^\d|v|^\d + q^\d|w|^\d.
 \]
Hence, if $|p-1|/q =\la$, we obtain the bound
\begin{equation}\nonumber
\displaystyle \frac d{dt}\int_{\real}|v|^{2+\delta}g(v,t)\, dv
+\frac{2+\delta}2q^2\left[ 1+\la^2 \right]
  \int_{\real}|v|^{2+\delta}g(v,t) \, dv \le
\end{equation}
\begin{equation}\nonumber
 + \frac 12 (2+\d)(1+\d)q^2 \int_{\real^2}
 \left((1+p)^\d|v|^\d + q^\d|w|^\d\right)(\la v+w)^2 g(v)
 g(w)\,dv\,dw,
\end{equation}
or, what is the same,
\begin{equation}\label{inv5}
\displaystyle \frac d{dt}\int_{\real}|v|^{2+\delta}g(v,t)\, dv \le
q^2 C(\la,q)
  \int_{\real}|v|^{2+\delta}g(v,t) \, dv \, .
\end{equation}
If we now use \fer{resc2}, it holds
\begin{equation}\label{inv6}
\displaystyle \frac d{d\tau}\int_{\real}|v|^{2+\delta}h(v,\tau)\,
dv \le  C(\la,q)
  \int_{\real}|v|^{2+\delta}h(v,\tau) \, dv \, ,
\end{equation}
namely the uniform boundedness of the $(2+\d)$-moment of
$h(v,\tau)$ with respect to $q$, for any fixed time $\tau$.

 Consider now the remainder \fer{resto}, which can be rewritten as
  \be\label{resto1}
 R(p, q) = \frac {q^2}{2}\int_{\real^2}\left(\frac{p-1}q v +w\right)^2  \left(
\phi''(\tilde v)- \phi''( v)\right)
 h(v)h(w)  dv\,d w .
 \ee
 We need the following

 \begin{definition}\label{reg}
 Let $\bigF_s(\real)$,  be the class of all real  functions $\phi$ on $\real$ such that
$\phi^{(m)}(v)$ is H\"older continuous of order $\delta$,
 \bq\label{lip} \|\phi^{(m)}\|_\delta= \sup_{v\not= w} \frac{|\phi^{(m)}(v) -\phi^{(m)}(w)|}{
|v-w|^\delta} <\infty,
 \eq
 the integer $m$ and the number $0 <\delta \le 1$ are such that $m+\delta =s$, and
$\phi^{(m)}$ denotes the $m$-th derivative of $g$.
\end{definition}

If $\phi \in \bigF_s(\real)$, with $s=2+\d$,
 \be
 \left| \phi''(\tilde v)- \phi''( v)\right| \le
 \|\phi''\|_\delta |\tilde v - v|^\d \le \|\phi''\|_\delta |(p-1)v
 +qw|^\d.
 \ee
 In this case,
 \be\nonumber
 R(p, q) \le \frac {q^{2+\d}}{2}\|\phi''\|_\delta \int_{\real^2}\left(\frac{p-1}q v
 +w\right)^{2+ \d}
 h(v)h(w)  dv\,d w \le
 \ee
 \be
\frac {q^{2+\d}}{2}\|\phi''\|_\delta
C_2(\la,q)\int_{\real^2}|v|^{2+ \d}h(v)\, dv.
 \ee
 Thanks to the uniform bound on $(2+\d)$-moment of
$h(v,\tau)$ , it follows that, for any fixed time $\tau >0$,
 \be
\lim_{q \to 0} \frac 1{q^2}\, R(p, q) = 0
 \ee
 as soon as  $\phi \in
\bigF_s(\real)$, with $s=2+\d$. This implies that the limit
equation is the Fokker-Planck equation \fer{FP}. We proved

\begin{theorem}\label{FP11}
Let the probability density $f_0 \in \bigM_\a$, where $\a= 2+
\delta$ for some $\delta>0$, and let the mixing parameters satisfy
 \[
 \frac{(p-1)^2}{q^2} = \la^2,
 \]
 for some constant $\la$ fixed.
Then, as $q \to 0$, for all $\phi \in \bigF_s(\real)$, with
$s=2+\d$ the weak solution to the Boltzmann equation \fer{quasiFP}
for the scaled density
 $h(v,\tau)=g(v,t)$, with $\tau = q^2 t$
converges, up to extraction of a subsequence, to a probability
density $h(w,\tau) $. This density  is a weak solution of the
Fokker-Planck equation \rf{FP}.
\end{theorem}

\subsection{A comparison of tails}\label{ta}

The result of Section \ref{grazz} establishes a rigorous
connection between the collisional kinetic equation \fer{eq:boltz}
and the Fokker--Planck equation \fer{FP}. The result of Lemma
\ref{SS}, coupled with the comment of Remark \ref{poss} then shows
that there exists a link between tails of the stationary solution
of Fokker--Plank and Boltzmann equations. In fact, one can choose
$\lambda^2>0$ in Theorem  \ref{FP1} if and only if the mixing
parameters $p$ and $q$ satisfy the conditions of the
aforementioned Lemma \ref{SS}. Since the reckoning of the size of
the tails is immediate in the Fokker--Planck case, it would be
important to know if one can extract from this knowledge
information about size of the tails of the Boltzmann equation.

Since the size of tails in the Boltzmann equation is given by the
positive root of the equation
 \[
\SS (\delta) = 0,
 \]
 where $\SS$ is the function \fer{key}, we will try to
extract information by comparing this root with the value of the
parameter $\lambda$ that characterizes the tails of the
Fokker--Planck equation. If $p>1$, using a Taylor expansion of
$\SS(\delta)$, with $p = 1+\la q$, we obtain
 \be\label{root}
 \frac{\SS(\delta)}{q^2} = \frac{2+\delta}2\left[ \left( \la^2\delta
 -1\right) + \frac 2{2+\delta}\, q^\delta + \frac{(1+\delta)\delta}3 \la^3 \frac{\bar
 q^3}{q^2} \right],
  \ee
where $ 0 \le \bar q \le q$. This shows that, in the scaling of
Theorem \ref{FP1}, the positive root $\delta^*(q)$ of $\SS(\delta)
= 0$ converges, as $q \to 0$ to the value $1/\la^2$, which
characterizes the tails of the Fokker--Planck equation. When $\la
>0$, one can easily argue that $\delta^*(q) < 1/\la^2$. In this case,
in fact,

\be\label{root2}
 \frac{\SS(\delta)}{q^2} = \frac{2+\delta}2\left[ \left( \la^2\delta
 -1\right) + A \right],
  \ee
where $A > 0$ if $q>0$. Hence
\be
 \frac{\SS(1/\la^2)}{q^2} = \frac{2+\delta}2 A  >0,
  \ee
that, by virtue of the convexity properties of $\SS(\delta)$
implies $\delta^*(q) < 1/\la^2$.

A weaker information can be extracted when $p<1$ while $p^2+q^2
<1$. In this case, writing $p = 1 -\la q$, $ \la >0$, we obtain
 \be\label{root3}
 \frac{\SS(\delta)}{q^2} = \frac{2+\delta}2\left[ \left( \la^2\delta
 -1\right) + \frac 2{2+\delta}\, q^\delta - \frac{(1+\delta)\delta}3 \la^3 \frac{\bar
 q^3}{q^2} \right],
  \ee
where $ 0 \le \bar q \le q$. Let us set
 \be\label{bbb}
 q \le \frac{B\la}{1+\la^2},
 \ee
 where $B \le 2$. In fact, when $p<1$ Lemma \ref{SS} implies that
 there is formation of tails only when $p$ and $q$ are such that
 $p^2+q^2 <1$, which is equivalent to the condition
 \be\label{bbb1}
 q < \frac{2\la}{1+\la^2}.
 \ee
Hence, when $q$ satisfies \fer{bbb}, from \fer{root3} we obtain
the inequality
 \be\label{root4}
 \frac{\SS(\delta)}{q^2} \ge \frac{2+\delta}2\left[ \left( \la^2\delta
 -1\right) - B\frac{(1+\delta)\delta}3  \frac{\la^4}{1+\la^2}
 \right].
  \ee
Easy computations then show that, if $\delta = r\la^2$, with
$0<r<1$, the right--hand side of \fer{root4} is nonnegative as
soon as
 \[
 3r(1-r)(1+\la^2) \ge B(1+r\la^2).
 \]
 Hence, the biggest value of $B$ for which the right--hand side of
 \fer{root4} is nonnegative is attained when $r =1/2$. In this
 case, $B= 3/4$, and $\delta^*(q) < 2/\la^2$. We can collect the
 previous analysis into the following

\begin{lemma}\label{tails11}
Let the mixing parameters satisfy
 \[
 \frac{(p-1)^2}{q^2} = \la^2,
 \]
for some constant $\la$ fixed. Then, if $p>1$ the positive root
$\delta^*(q)$ of the equation $\SS (\delta) = 0$, characterizing
the tails of the Boltzmann equation, satisfies the bound
$\delta^*(q) < 1/\la^2$. If $p>1$, and at the same time $q$
satisfies the bound \fer{bbb} with $B= 3/4$, the positive root
$\delta^*(q)$ of the equation $\SS (\delta) = 0$, satisfies the
bound $\delta^*(q) < 2/\la^2$.
\end{lemma}

We remark here that, in the case $p <1$, setting $\delta = 1$ we
obtain an exact formula for $\SS(1)$,
 \be
 \frac{\SS(1)}{q^2} = \frac{3}2\left[ \left( \la^2
 -1\right) + \frac 2{3}\, q - \frac{2}3 \la^3 q \right].
  \ee
Choosing $\la =1$, we get $\SS(1)=0$. This case, that corresponds
to the conservation of momentum in the Boltzmann equation has
tails which are invariant with respect to $q$ (see Remark
\ref{granu}).


\subsection{Kinetic models of economy}

The analysis of Sections \ref{appl}, \ref{uni}, \ref{convv} and
\ref{grazz} can be easily extended to equation \fer{eq:bolt1} for
the wealth distribution. We can in fact resort to the methods
introduced for the kinetic equation on the whole real line simply
setting
 \be
 F(v,t) = f(v,t)I(v \ge 0), \quad v \in \real,
 \ee
 where $I(A)$ is the indicator function of the set $A$.
 With this notation, equation \fer{eq:bolt1} can be rewritten as
 equation \fer{eq:boltz},
 \be
 \frac{\partial F(v)}{\partial t} =
\int_{\real}\left( \frac 1J F(v_*)
 F(w_*) - F(v) F(w)\right) d
 w. \label{eq:bolt11}
  \ee
 Likewise, the weak form \fer{eco} reads
  \be\label{eco1}
\frac d{dt}\int_{\real} \phi(v)F(v,t)\,dv   =
 \int_{\real^2} F(v,t)F(w,t) ( \phi(v^*)-\phi(v)) dv d w .
 \ee
We recall that the role of the energy is now supplied by the mean
$m(t) = \int vF(v,t)\, dv$. To look for self--similarity we scale
our solution according to
 \be\label{sca1}
  G(v,t) = { m(t)}F\left( { m(t)}v, t \right),
  \ee
which implies that $\int v G(v,t)= 1$ for all $t \ge 0$. Hence,
without loss of generality, if we fix the initial density to
satisfy
 \be\label{norm2}
 \int_{\real} F_0(v)\, dv =1 \, ; \quad    \int_{\real} v F_0(v)\, dv =  1
 \, ,
 \ee
the solution $G(v,t)$ satisfies \fer{norm2}. Then, the same
computations of Section \ref{appl} show the following

\begin{theorem} \label{contr1}
Let $f_1(t)$ and $f_2(t)$ be two solutions of the Boltzmann
equation~\fer{eq:bolt1}, corresponding to  initial values
$f_{1,0}$ and $f_{2,0}$ satisfying conditions \fer{norm2}. Then,
if for some  $1 \le s \le 2$,  $d_s(f_{1,0},f_{2,0})$ is
 bounded, for all times $t \geq 0$,
 \be\label{dec31}
d_s(f_1(t), f_2(t)) \leq \exp\left\{ (p^s + q^s
-1)t\right\}d_s(f_{1,0},f_{2,0}).
  \ee
In particular, let $f_0$ be a nonnegative density satisfying
conditions \fer{norm1}. Then, there exists a unique weak solution
$f(t)$ of the Boltzmann equation, such that $f(0) = f_0$. In case
$p^s + q^s -1 <0$ the distance $d_s$ is contracting exponentially
in time.
\end{theorem}

 Since by \fer{sca1}
 \[
 \widehat{G}(\xi) = \widehat{G}\left(\frac\xi{
 m(t)}\right),
 \]
from \fer{ds1} we obtain the bound
 \begin{equation} \label{ds12}
 d_s (g_1(t),g_2(t))= \sup_{\xi \in {\real}} \frac{|\widehat{G_1}(\xi,t) -
\widehat{G_2}(\xi,t)|}{|\xi|^s} =  \left(\frac 1{
 m(t)}\right)^s d_s (f_1(t),f_2(t)).
 \end{equation}
 Using \fer{dec31}, we finally conclude that, if $g_1(t)$ and
 $g_2(t)$ are two solutions of the scaled Boltzmann
equation~\fer{eq:bolt1}, corresponding to  initial values
$f_{1,0}$ and $f_{2,0}$ satisfying conditions \fer{norm2},
 Then, if $1 \le s \le 2$,
for all times $t \geq 0$,
 \be\label{dec11}
  d_s(g_1(t), g_2(t)) \leq \exp\left\{\left[ (p^s + q^s -1)
   -s(p+q -1)\right]t\right\}d_s(f_{1,0},f_{2,0}).
  \ee
Let us define, for $\d \ge 0$,
 \be\label{key1}
 \RR (\delta) = p^{1+\delta} + q^{1+\delta} -1 -(1+\d)\left(p+q
 -1\right).
 \ee
Then, the sign of $\RR$ now determines the asymptotic behavior of
the distance $d_s(g_1(t), g_2(t))$. With few differences, the
proof leading to Lemma \ref{SS} can be repeated, obtaining

\begin{lemma}\label{RR}
 Let $\RR(\d), \d \ge 0$ be the function defined by \fer{key1}.
 Given a constant $\la >0$, if $p+q < 1$, let us define $p=1- \la\sqrt q$ . Then,
 provided $q < 1/\la^2$ there exists an interval $I_-= (0 , \bar \d_-(q))$
 such that $\RR(\d)<0$ for $\d \in I_-$. If $p+q > 1$, and $p=1+ \la
 \sqrt q$ there exists a interval $I_+= (0 , \bar \d_+(q))$
 such that $\RR(\d)<0$ for $\d \in I_+$. In the remaining cases, namely when $p+q =
 1$ or $p+q > 1$ but $p<1$, $\RR(\d)<0$ for all $\d>0$.
 \end{lemma}

The main consequence of Lemma \ref{RR} is contained into the
following.

\begin{theorem} \label{de1}
Let $g_1(t)$ and $g_2(t)$ be two solutions of the Boltzmann
equation~\fer{eq:bolt1}, corresponding to  initial values
$f_{1,0}$ and $f_{2,0}$ satisfying conditions \fer{norm2}. Then,
there exists a constant $\bar\d >0$ such that, if $1 < s < 1+
\bar\d$, for all times $t \geq 0$,
 \be\label{decc1}
  d_s(g_1(t), g_2(t)) \leq \exp\left\{ -C_st\right\}d_s(f_{1,0},f_{2,0}).
  \ee
 The constant $C_s = -\RR(s-1)$ is strictly positive, and
 the distance $d_s$ is contracting exponentially
in time.
\end{theorem}

 Existence and uniqueness of the stationary solution to equation
 \fer{weak boltz11} follows along the same lines of Section
 \ref{convv}.
 The main result is now contained into the following.

 \begin{theorem} \label{pnot11}
Let $\d >0$ be such that $\RR(\d)<0$, and let $g_{\infty}(v)$ be
the unique stationary solution to equation \fer{weak boltz11}. Let
$g(v,t)$ be the weak solution of the Boltzmann equation~\fer{weak
boltz11}, corresponding to the initial density $f_0$ satisfying
\[
\int_{\real_+} |v|^{1+\delta} \, f_0(v)\, dv<\infty .
\]
 Then, $g(v,t)$ satisfies
\[
\int_{\real_+} |v|^{1+\delta} \, g(v, t)\, dv \le c_{\delta}
<\infty ,
\]
for some constant $c_\d$ depending only on $p$ and $q$.
 Moreover,  $g(v,t)$ converges exponentially fast in
Fourier metric towards $g_{\infty}(v)$, and the following bound
holds
\begin{equation}\label{conv11}
d_{1+\delta}(g(t),g_{\infty} ) \le d_{1+\delta}(f_0,g_{\infty}
)\exp\left\{-|\RR(\delta)|t\right\}
\end{equation}
where $\RR(\delta)$ is given by \fer{key1}.
\end{theorem}

Depending on the values of the mixing parameters $p$ and $q$, the
stationary solution $g_{\infty}$ can have overpopulated tails. The
Fourier transform of $g_{\infty}$ satisfies the equation
\begin{equation}\label{fougg}
-(p+q-1)\xi \frac{\partial{\widehat{G}}}{\partial \xi}+\widehat
G(\xi)=\widehat G(p\xi)\widehat G(q\xi).
\end{equation}
We set
 \be\label{serie1}
\widehat{G}(\xi) = 1-i\xi + A|\xi|^{1+\delta} + \dots
 \ee
 which takes into account the fact that  $g_{\infty}$ satisfies
 conditions \fer{norm2}. The leading small $\xi$-behavior of the
 singular component will reflect an algebraic tail of the velocity
 distribution. Substitution of expression \fer{serie1} into
 \fer{fougg} shows that the coefficient of the power
 $|\xi|^{1+\delta}$ is $A\RR(\d)$. Thus, the term
 $A|\xi|^{1+\delta}$  can appear in the expansion of
 $\widehat{G}(\xi)$ as soon as $\d$ is such that $\RR(\d) = 0, \d >0$.
As before, tails in the stationary distributions are present in
all cases in which there exists a $\d = \bar \d>0$ such that
$\RR(\bar\d)=0$. Now the answer is contained into Lemma \ref{RR}.

Last, one can justify rigorously the passage to the Fokker-Planck
equation \fer{FP2b}.

\begin{theorem}
Let the probability density $f_0 \in \bigM_\a$, where $\a= 1+
\delta$ for some $\delta>0$, and let the mixing parameters satisfy
 \[
 \frac{(p-1)^2}{q} = \la,
 \]
 for some $\la >0$ fixed.
Then, as $q \to 0$, for all $\phi \in \bigF_s(\real)$, with
$s=1+\d$ the weak solution to the Boltzmann equation \fer{quasiFF}
for the scaled density
 $h(v,\tau)=g(v,t)$, with $\tau = q t$
converges, up to extraction of a subsequence, to a probability
density $h(w,\tau) $. This density  is a weak solution of the
Fokker-Planck equation \rf{FP2b}.
\end{theorem}

 We finally remark that the discussion of Section \ref{ta}, with minor
 modifications, can be adapted to establish connections between
 the size of the tails of the kinetic and Fokker--Planck models.


\section{Numerical examples} In this paragraph, we shall
compare the self--similar stationary results obtained by using
Monte Carlo simulation of the kinetic model with the stationary
state of the Fokker-Planck model. The method we adopted is based
on Bird's time counter approach at each time step followed by a
renormalization procedure according to the self-similar scaling
used. We refer to \cite{P2} for more details on the use of Monte
Carlo method for Boltzmann equations.


We used $N=5000$ particles and perform several iterations until a
stationary state is reached. The distribution is then averaged
over the next $4000$ iterations in order to reduce statistical
fluctuations. Clearly, due to the slow convergence of the Monte
Carlo method near the tails, some small fluctuations are still
present for large velocities.

\subsubsection*{Gaussian behavior} 
\begin{figure}[htb]
\begin{center}
\includegraphics[scale=.42]{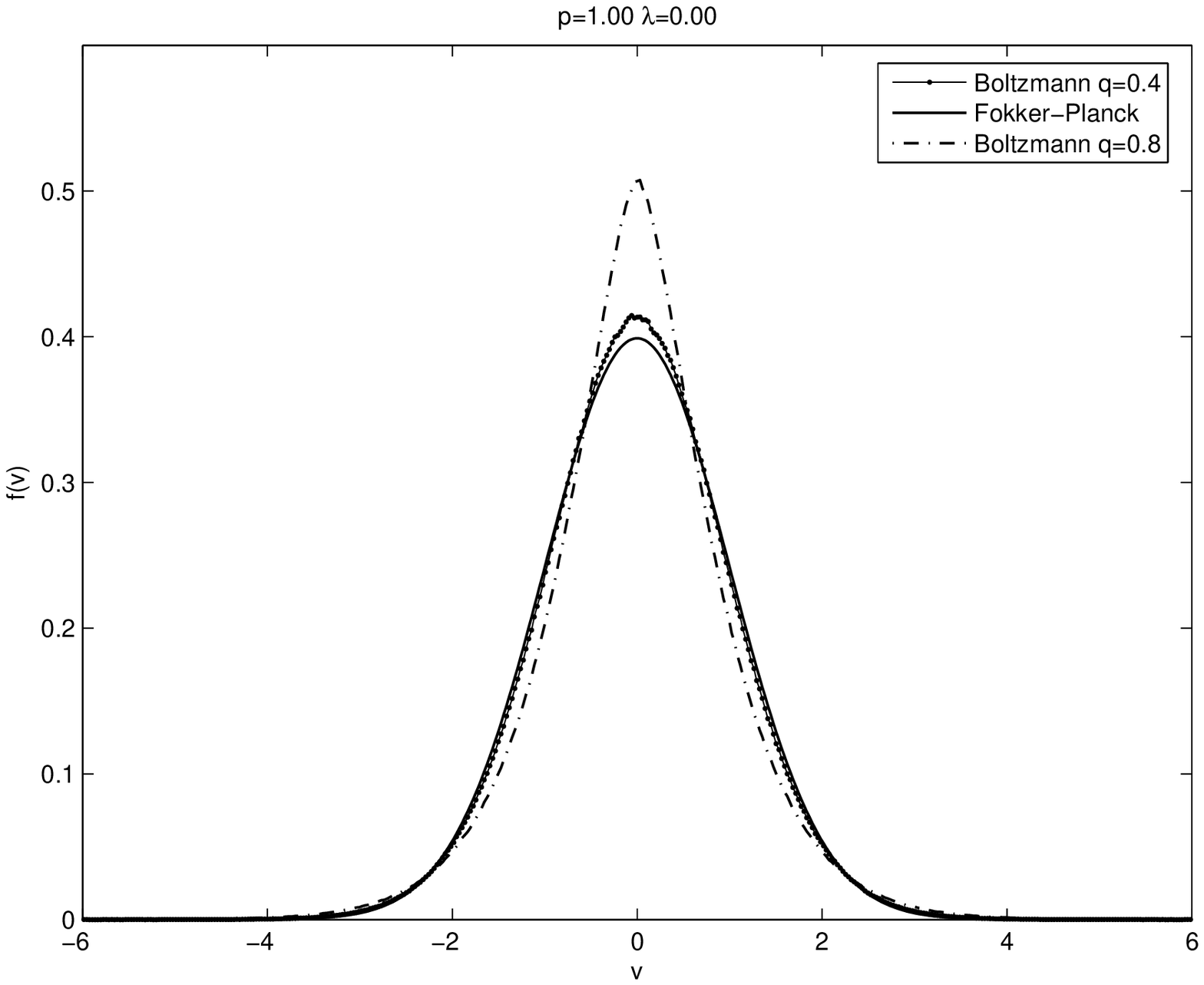}
\includegraphics[scale=.42]{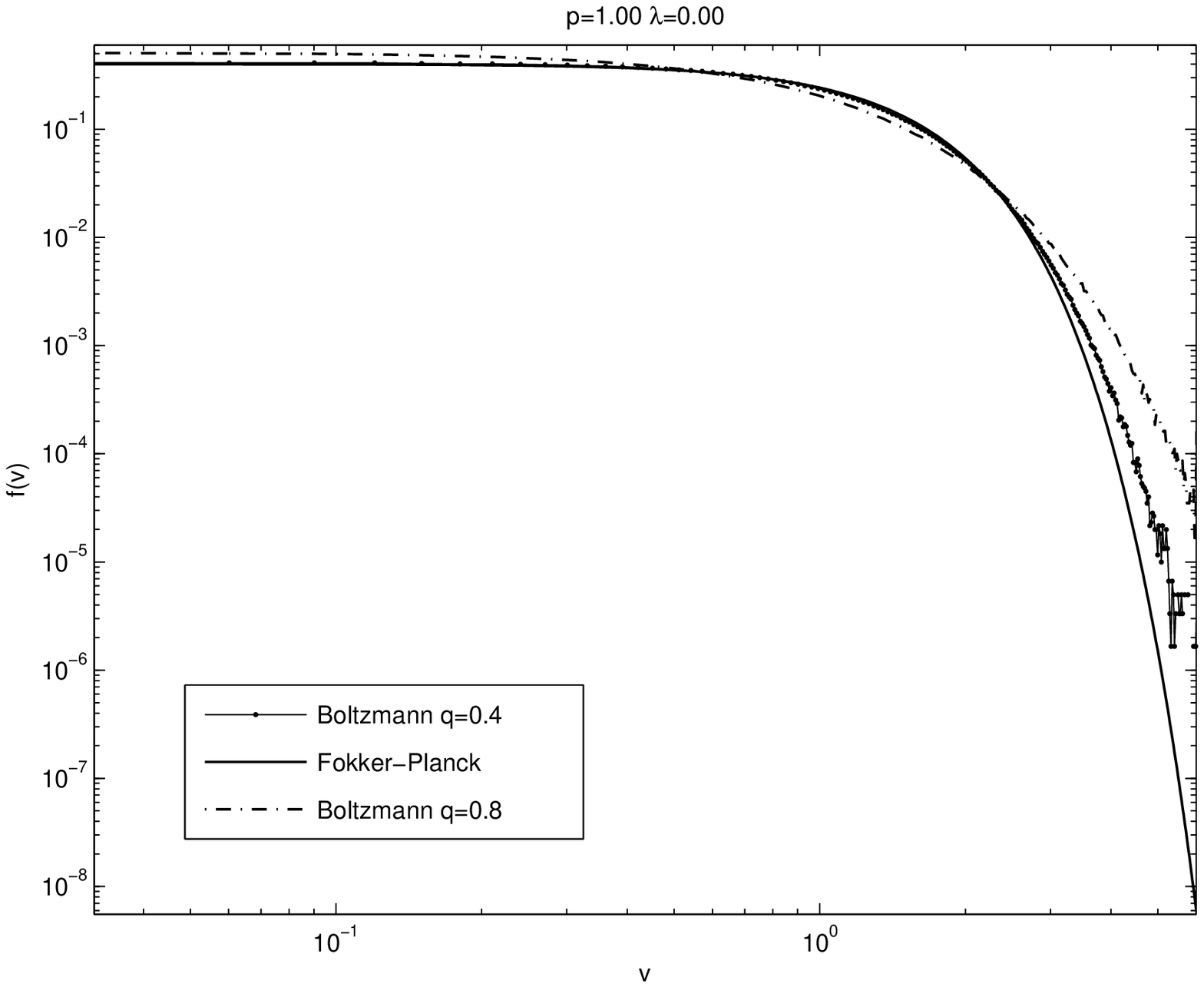}
\end{center}
\caption{Asymptotic behavior for $\lambda=0$ of the Fokker-Planck
model and the Boltzmann model with $p=1$ and $q=0.4, 0.8$. Figure
on the right is in loglog-scale.} \label{fig:11}
\end{figure}

First we consider the case
$\lambda=0$ for which the steady state of the Fokker-Planck
asymptotic is the Gaussian (\ref{eq:max}). We fix $p=1$ so that
for $q < 1/\sqrt{2}$ we expect Gaussian behavior also in the
kinetic model. We report the results obtained for $q=0.4$ and
$q=0.8$ in Figure \ref{fig:11}.

\subsubsection*{Formation of power laws} Next we simulate the
formation of power laws for positive $\lambda$. We take $p=1.2$
and $q=.4$ which correspond to $\lambda=0.5$. Keeping the same
value of $\lambda$ we then take $q=0.1$ and $p=1.05$. In Figure
\ref{fig:11} we plot the results showing convergence towards the
Fokker-Planck behavior.

\begin{figure}[htb]
\begin{center}
\includegraphics[scale=.42]{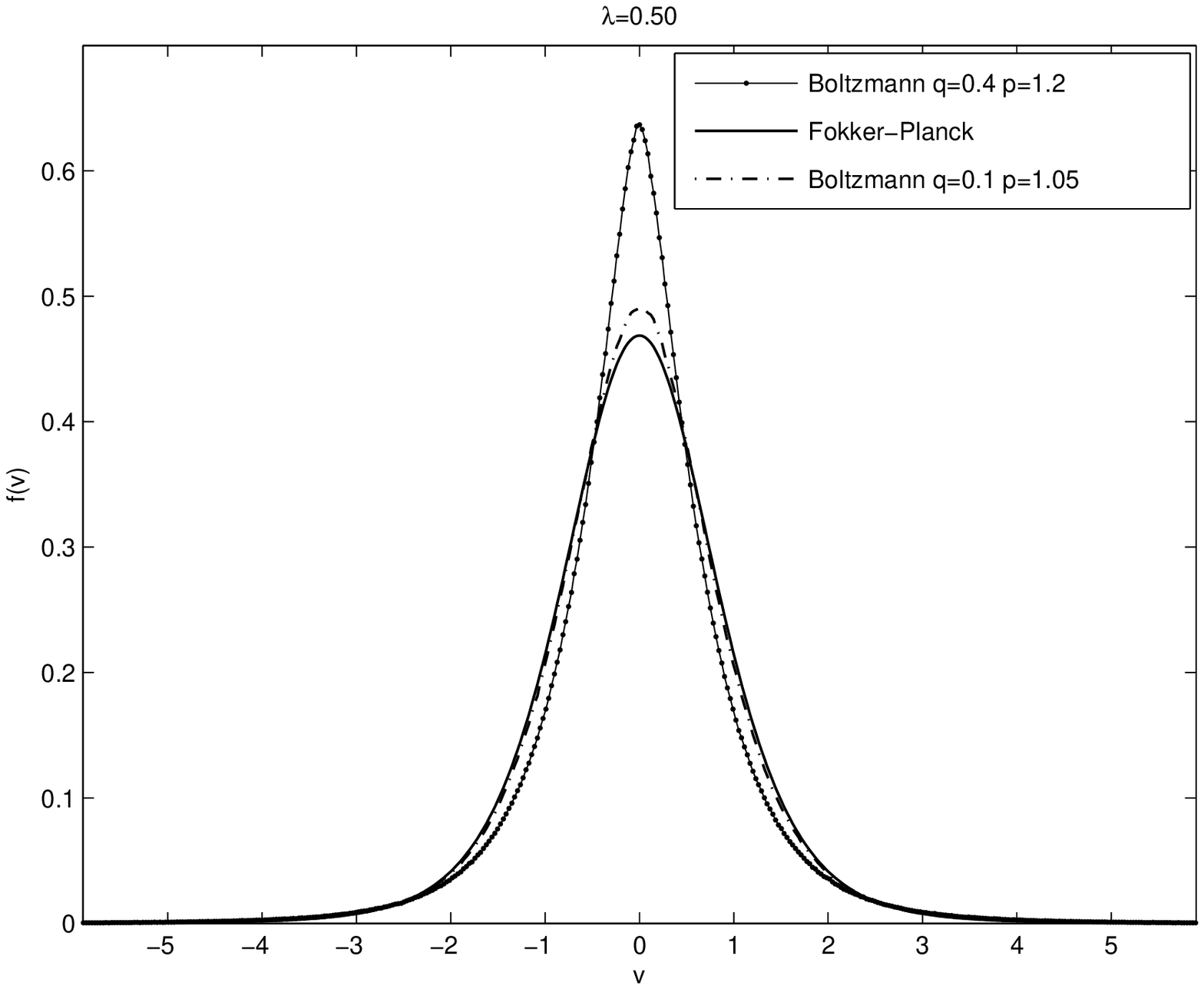}
\includegraphics[scale=.42]{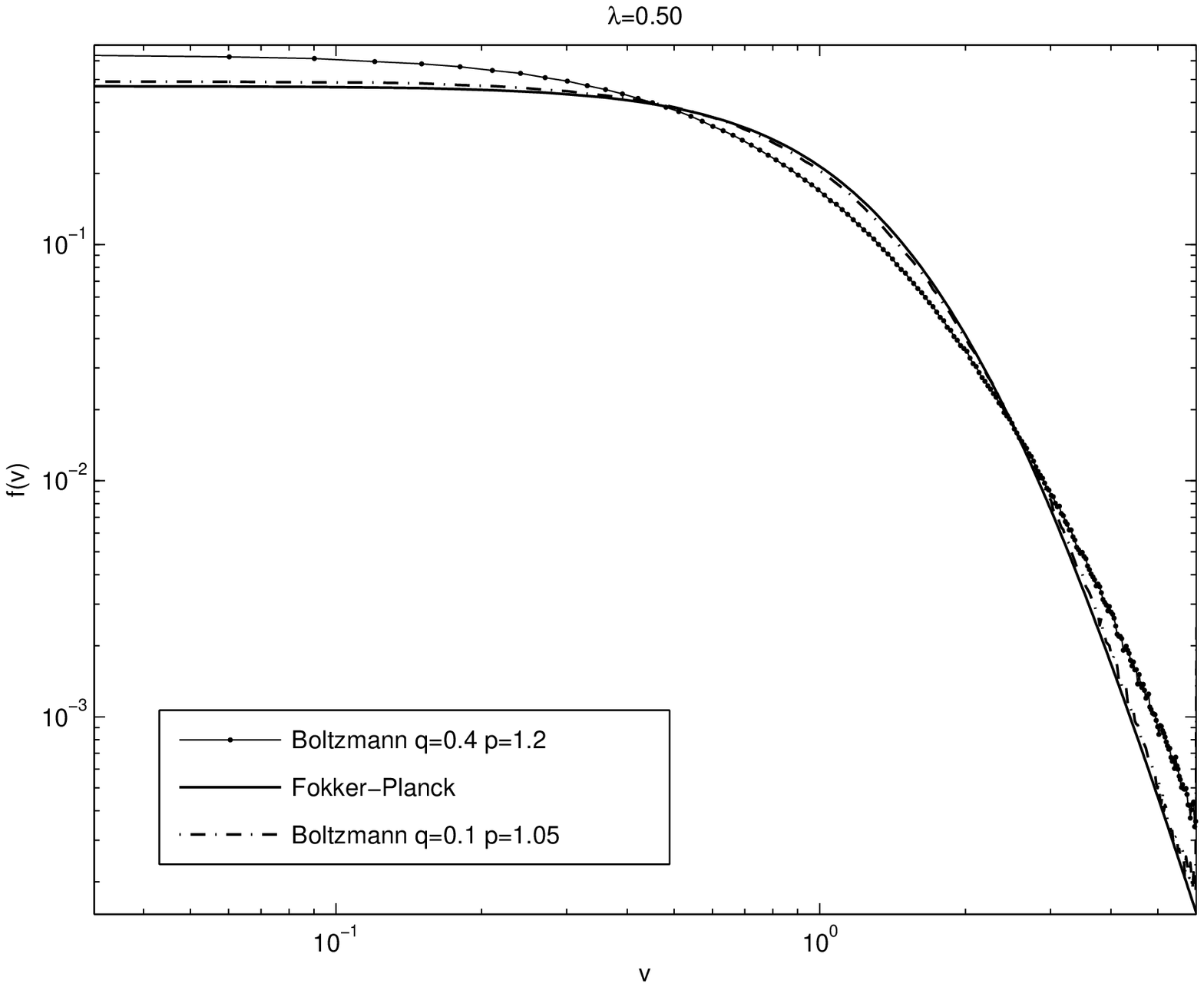}
\end{center}
\caption{Asymptotic behavior for $\lambda=0.5$ of the
Fokker-Planck model and the Boltzmann model for $p=1.2$, $q=0.4$
and $p=1.05$, $q=0.1$. Figure on the right is in loglog-scale.}
\label{fig:12}
\end{figure}

\subsubsection*{A simple growing economy}

\begin{figure}[ht]
\begin{center}
\includegraphics[scale=.42]{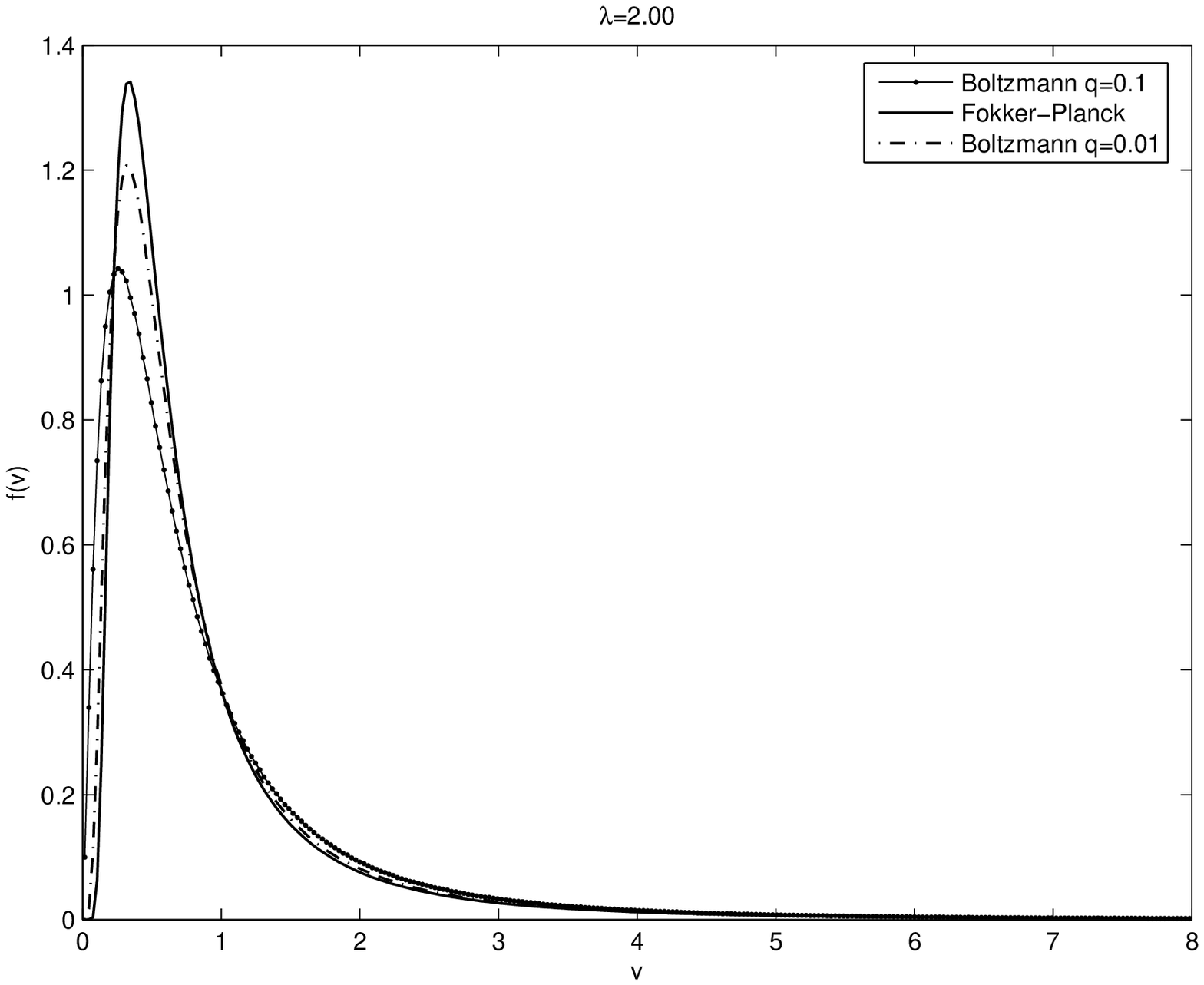}
\includegraphics[scale=.42]{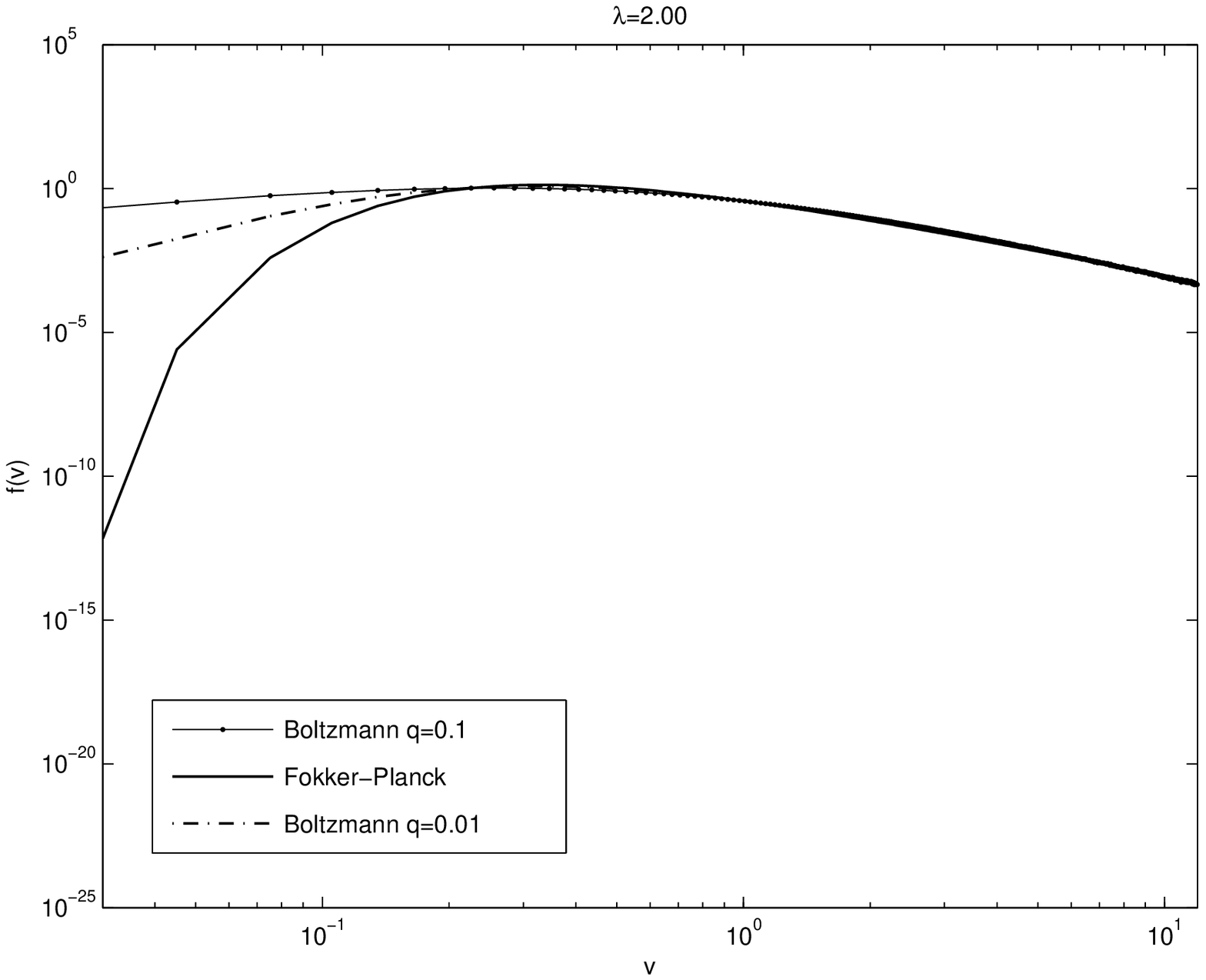}
\end{center}
\caption{Asymptotic behavior for $\lambda=2$ of the Fokker-Planck
model and the Boltzmann model for $p=1-q+2\sqrt{q}$, $q=0.1$ and
$q=0.01$. Figure on the right is in loglog-scale.} \label{fig:2}
\end{figure}

We take the case of a growing economy for $p=1-q+2\sqrt{q}$ thus
corresponding to the limit Fokker-Planck steady state (\ref{equi})
with $\lambda=2$ and $\mu=2$. As prescribed from our theoretical
analysis we observe that the equilibrium distribution converges
toward the Fokker-Planck limit as $q$ goes to 0, with $\lambda$
fixed. The results are reported in Figures \ref{fig:2}.


\section{Conclusions}

In this paper we studied the large--time behavior of a simple
one-dimensional kinetic model of Maxwell type, in two situations,
depending wether the velocity variable can take values on $\real$
or in $\real_+$, the former case describing nonconservative models
of kinetic theory of rarefied gases, the latter elementary kinetic
models of open economies. In both situations it has been shown
that the lack of conservation laws leads to situations in which
the self--similar solution has overpopulated tails. This is
particularly important in the case of economy, where elementary
explanations of the formation of Pareto tails can help to handle
more complex models of society wealth distribution, where various
other factors occur. It would be certainly interesting to extend a
similar analysis to more realistic situations. Recently, a kinetic
model including market returns has been introduced \cite{CPT}.
While for this model the asymptotic convergence to the
Fokker--Planck limit can be obtained, the property of creation of
overpopulated tails has been shown only by numerical simulation.
In realistic models, in fact, there is a strong correlation among
densities, due to the constraint of having non-negative wealths
after trades, and this appears difficult to treat from a
mathematical point of view. A further point deserves to be
mentioned. Recent studies have shown that, while overpopulated
tails seem to be generic feature of the non-conservative collision
mechanism, in the kinetic theory of the Boltzmann equation
power--like tails only occur in the {\emph{ borderline}} case of
Maxwell molecules interactions \cite{BC4, BCT, EB1, EB2}, whereas
in general collision dissipative processes have stretched
exponential tail behaviors \cite{BGP}. It could be conjectured
that the corresponding phenomenon in general kinetic models of
economy with wealth--depending collision frequency manifest a
behavior in the form of a lognormal type distribution.


\bigskip \noindent {\bf Acknowledgment:} The authors acknowledge support from the IHP
network HYKE ``Hyperbolic and Kinetic Equations: Asymptotics, Numerics, Applications''
HPRN-CT-2002-00282 funded by the EC., and from the Italian MIUR, project ``Mathematical Problems
of Kinetic Theories''.

%

%

\begin{thebibliography}{[KLR73]}

\bibitem[BMP02]{BMP}
Baldassarri A., Marini Bettolo Marconi U.,  Puglisi A. :
\newblock Kinetic models of inelastic gases.
\newblock {\em Mat. Mod. Meth. Appl. Sci.} \textbf{12}  965--983 (2002).

\bibitem[BBLR03]{BBLR}
Ben-Avraham D., Ben-Naim E.,  Lindenberg K.,  Rosas A.:
\newblock Self-similarity in random collision processes,
\newblock  {\em Phys. Rev. E}, \textbf{68}, R050103 (2003).

\bibitem[BK00]{BK}
Ben-Naim E., Krapivski P. :
\newblock Multiscaling in inelastic collisions.
\newblock {\em Phys. Rev. E}, \textbf{61}, R5--R8 (2000).

\bibitem[Bo88]{Bob}
Bobylev A.V.:
\newblock The theory of the nonlinear spatially uniform Boltzmann equation for Maxwellian
molecules. {\em Sov. Sci. Rev. c} \textbf{7}, (1988) 111-233.


\bibitem[BCG00]{BCG}
Bobylev A.V., Carrillo J.A., Gamba I.:
\newblock On some properties of kinetic and hydrodynamics
equations for inelastic interactions.
\newblock {\em  J. Statist.
Phys.} \textbf {98},  743--773 (2000).



\bibitem[BC03]{BC4}
Bobylev A.V., Cercignani C.:
\newblock Self-similar asymptotics for the Boltzmann equation with inelastic and elastic interactions.
\newblock {\em  J. Statist.
Phys.} \textbf {110}, 333-375 (2003).


\bibitem[BCT03]{BCT}
 Bobylev A.V., Cercignani C.,  Toscani G. :
\newblock Proof of an asymptotic property of self-similar solutions of the Boltzmann equation for granular materials.
\newblock {\em  J. Statist.
Phys.} \textbf {111},   403-417 (2003).

\bibitem[BGP04]{BGP}
 Bobylev A.V., Gamba I.,  Panferov V.A. :
\newblock Moment inequalities and high--energy tails for Boltzmann equations with inelastic interactions.
\newblock {\em  J. Statist.
Phys.} \textbf {116},   1651-1682 (2004).

\bibitem[BM00]{BM}
Bouchaud J.P., M\'ezard M.:
\newblock{Wealth condensation in a simple model of economy}, {\em Physica A},
\textbf{282}, 536- (2000).


\bibitem[CGT99]{CGT}
Carlen E.A., Gabetta E., Toscani .:
\newblock  Propagation of smoothness and the rate of exponential convergence to
equilibrium for a spatially homogeneous Maxwellian gas.
\newblock {\em Commun. Math. Phys.}{\bf 305}, 521-546 (1999).

\bibitem[CCG00]{CCG1}
Carlen E.A., Carvalho M.C., Gabetta E.:
\newblock  Central limit theorem for Maxwellian molecules and truncation of the
Wild expansion.
\newblock {\em Commun. Pure Appl. Math.},  \textbf{ 53} (2000), 370--397.

\bibitem[CIP94]{CIP}
C. Cercignani, R. Illner, M. Pulvirenti,
\newblock {\it The mathematical theory of dilute gases}.
\newblock {\em Springer Series in Applied Mathematical Sciences},
Vol.\textbf{  106} Springer--Verlag, New York 1994.


\bibitem[CPT04]{CPT}
Cordier S.,  Pareschi L., Toscani G.:
\newblock On a kinetic model for a simple market economy.
\newblock {\em J. Stat. Phys.}, to appear, (2005).

\bibitem[EB02a]{EB1}
 Ernst M.H.,   Brito R.:
\newblock High energy tails for inelastic Maxwell models.
{\em Europhys. Lett}, \textbf{ 43}, 497-502 (2002).

\bibitem[EB02b]{EB2}  Ernst M.H.,   Brito R.:
\newblock Scaling solutions of inelastic Boltzmann equation with over-populated high energy tails.
{\em J. Statist. Phys.}, \textbf{109},  407-432 (2002).


\bibitem[DY00]{YD}
Dr\v{a}gulescu A., Yakovenko V.M.:
\newblock{Statistical mechanics of money}, {\em Eur. Phys. J. B.}, \textbf{17}, 723--729
(2000).

\bibitem[GGPS03]{GGPS}
Gabaix X., Gopikrishnan P.,  Plerou V., Stanley H.E.:
\newblock{A Theory of Power-Law Distributions in Financial Market
Fluctuations}, {\em Nature} \textbf{423}, 267--270 (2003).

\bibitem[GTW95]{GTW}
Gabetta E., Toscani G., Wennberg, B.:
\newblock  Metrics for probability distributions and the trend to equilibrium
for solutions of the Boltzmann equation,
\newblock {\em J. Stat. Phys.}, { \bf 81 } (1995) {901--934}.

\bibitem[GJT02]{GJT}
Goudon T., Junca S., Toscani G.:
\newblock Fourier-based distances and Berry-Esseen like
inequalities for smooth densities. {\em Monatsh. Math.}, {\bf 135
}(2002)  115-136.

\bibitem[IKR98]{IKR}
Ispolatov S., Krapivsky P.L., Redner S.:
\newblock{Wealth distributions in asset exchange models}, {\em Eur. Phys. J. B}, 2,
267--276 (1998).

\bibitem[Ka59]{Kac}
Kac M.:
\newblock {\it Probability and related topics in the physical
sciences}, Interscience Publishers, London-New York 1959.

\bibitem[LR79]{LR}
Laha R.G., Rohatgi V.K.:
\newblock {\it Probability Theory},
John Wiley and Sons, New York, 1979.


\bibitem[BMRS02]{BMRS}Malcai O., Biham O.,
Solomon S., Richmond P. :
\newblock Theoretical analysis and simulations of the generalized Lotka-Volterra model, {\em Phys.
Rev. E}, \textbf{66}, 031102 (2002).

\bibitem[MK66]{McK}
 McKean, H.P. Jr.:
\newblock Speed of approach to equilibrium for Kac's caricature of a
Maxwellian gas.
\newblock {\em Arch. Rat. Mech. Anal.}, { \bf 21}:{343--367}, 1966.


\bibitem[MY93]{McNY}
 McNamara S.,  Young W.R.:
\newblock Kinetics of a one--dimensional granular
  medium in the quasi--elastic limit,
\newblock {\em Phys. Fluids A} \textbf{ 5},  34--45  (1993).

\bibitem[P04]{P}
 Pareschi L.:
\newblock Microscopic dynamics and mesoscopic modelling of a market economy, preprint
2004.

\bibitem[PR01]{P2}
Pareschi L., Russo G.: An introduction to Monte Carlo methods for
the Boltzmann equation. {\em ESAIM: Proceedings} {\bf 10}, 35--75,
(2001).


\bibitem[P897]{Par}
Pareto V.:
\newblock {\em Cours d'Economie Politique}, Lausanne and Paris, 1897.

\bibitem[PGS03]{PGS}
 Plerou V.,  Gopikrishnan P., Stanley H.E.:
\newblock {Two-Phase Behaviour of Financial Markets}, {\em Nature}
\textbf{421}, 130 (2003).

\bibitem[PB03]{BP}
Potters M., Bouchaud J.P.:
\newblock{Two--Phase behavior of financial markets}
arXiv: cond--mat/0304514 (2003)

\bibitem[PT03]{PT}
 Pulvirenti A., Toscani G. :
\newblock Asymptotic properties of the inelastic Kac model
\newblock {\em  J. Statist. Phys.} \textbf{ 114},  1453-1480 (2004).

\bibitem[S03]{Sl}
Slanina F.:
\newblock{Inelastically scattering particles and wealth distribution in an open economy}, preprint (2003)
cond-mat/0311025.

\bibitem[So98]{So}
Solomon S.: \newblock{Stochastic Lotka-Volterra systems of competing auto-catalytic agents lead generically to
truncated Pareto power wealth distribution, truncated Levy distribution of market returns, clustered volatility,
booms and crashes}, {\em Computational Finance 97}, eds. A-P. N. Refenes, A.N. Burgess, J.E. Moody (Kluwer Academic
Publishers 1998).

\bibitem[T00]{To}
 Toscani G. :
\newblock One-dimensional kinetic models of granular flows,
\newblock {\em RAIRO Mod\'el. Math. Anal. Num\'er.}
\textbf{  34},  1277-1292  (2000).

\bibitem[TV99]{TV1}
 Toscani G., Villani C.:
 \newblock Probability metrics and uniqueness of the solution to the Boltzmann
 equation for a Maxwell gas. {\em J. Statist. Phys.}, {\bf 94} (1999) 619-637.

\bibitem[V98a]{Vi}
 Villani C. :
\newblock  Contribution {\`a} l'{\'e}tude math{\'e}matique des
  {\'e}quations de {B}oltzmann et de {L}andau en th{\'e}orie cin{\'e}tique des
  gaz et des plasmas.
\newblock {\em PhD thesis, Univ. Paris-Dauphine}, (1998).

\bibitem[V98b]{vil2}
 Villani C.:
\newblock{On a new class of weak solutions to the spatially homogeneous Boltzmann
and Landau equations}, \newblock{\em Arch. Rational Mech. Anal}. {\bf 143},  273--307 (1998).

\end{thebibliography}
%


\end{document}